\documentclass[12pt]{amsart}
 \pdfoutput=1

\usepackage{amssymb, amsmath, amsthm, amsxtra, amscd}

\usepackage{multirow} 
\usepackage{graphicx}
\usepackage{pdfsync}
\usepackage{tikz}
\usetikzlibrary{cd}

\usepackage[linesnumbered]{algorithm2e}

\graphicspath{ {./pictures/} }

\newcommand{\calD}{\ensuremath{\mathcal{D}} }

\newcommand{\calF}{\ensuremath{\mathcal{F}} }

\newcommand{\calH}{\ensuremath{\mathcal{H}} }

\newcommand{\calL}{\ensuremath{\mathcal{L}} }

\newcommand{\calO}{\ensuremath{\mathcal{O}} }
\newcommand{\calP}{\ensuremath{\mathcal{P}} }
\newcommand{\calQ}{\ensuremath{\mathcal{Q}} }
\newcommand{\calR}{\ensuremath{\mathcal{R}} }
\newcommand{\calS}{\ensuremath{\mathcal{S}} }

\newcommand{\calT}{\ensuremath{\mathcal{T}} }

\newcommand{\la}{\langle}
\newcommand{\ra}{\rangle}

\newcommand{\bfa}{\ensuremath{\mathbf a} }
\newcommand{\bfb}{\ensuremath{\mathbf b} }

\newcommand{\bfe}{\ensuremath{\mathbf e} }

\newcommand{\bfn}{\mathbf n}

\newcommand{\bfp}{\mathbf{p}}
\newcommand{\bfq}{\mathbf{q}}

\newcommand{\bfsigma}{\ensuremath{\, \mathbf \sigma} }

\newcommand{\bfw}{\ensuremath{\mathbf w} }
\newcommand{\bfx}{\mathbf x}
\newcommand{\bfy}{\ensuremath{\mathbf y} }

\newcommand{\bfzero}{\mathbf 0}

\newcommand\boldP{\mathbb P}
\newcommand\boldR{\mathbb R}
\newcommand\boldN{\mathbb N}

\newcommand{\Dom}{\operatorname{Dom}}

\newcommand{\Vor}{\operatorname{Vor}}

\newcommand{\dist}{\operatorname{dist}}

\newcommand{\Cell}{\operatorname{Cell}}
\newcommand{\MD}{\operatorname{MD}}
\newcommand{\MRN}{\operatorname{MR^N}}
\newcommand{\MRL}{\operatorname{MR^\calL}}
\newcommand{\MRF}{\operatorname{MR^\calF}}
\newcommand{\MRI}{\operatorname{MR(I)}}
\newcommand{\Pow}{\operatorname{pow}}

\newcommand{\prc}{\operatorname{pr_c}}
\newcommand{\setint}{\operatorname{int}}
\newcommand{\setext}{\operatorname{ext}}
\newcommand{\Sph}{\operatorname{Sph}}

\newcommand{\D}{\displaystyle}

\newcommand{\smallfrac}[2]{\textstyle{\frac{#1}{#2}}  }

\usepackage{colortbl}

\theoremstyle{plain}
\newtheorem{thm}{Theorem} 
\newtheorem{prop}[thm]{Proposition}

\newtheorem{cor}[thm]{Corollary}

\theoremstyle{definition} 
\newtheorem{defn}[thm]{Definition}
\newtheorem{example}[thm]{Example}

\theoremstyle{remark}
\newtheorem{remark}[thm]{Remark}

\newtheorem*{conventions}{Conventions}
\newtheorem*{acknowledgments}{Acknowledgments}

\numberwithin{equation}{section}
\numberwithin{thm}{section}

\begin{document}

\title{Generalized Voronoi Diagrams and Lie Sphere Geometry}

\author{John Edwards}
\address{Department of Computer Science, Utah State, 4205 Old Main Hall
Logan, UT 84322-4205}
\email{john.edwards@usu.edu}

\author{Tracy Payne}
\address{Department of Mathematics and Statistics, Idaho State
  University, 921 S.\  8th St., 
Pocatello, ID 83209-8085}
\email{tracypayne@isu.edu}

\author{Elena Schafer}
\address{Department of Mathematics and Statistics,
  Idaho State University, 921 S.\  8th St., 
Pocatello, ID 83209-8085}
\email{elenaschafer12@gmail.com}

\begin{abstract}
We use Lie sphere geometry to describe two large categories of
generalized Voronoi diagrams that can be encoded in terms of the Lie
quadric, the Lie inner product, and polyhedra.  The first class
consists of diagrams defined in terms of extremal spheres in the space
of Lie spheres, and the second class includes minimization diagrams for
functions that can be expressed in terms of affine functions on a
higher-dimensional space.  These results unify and generalize previous
descriptions of generalized Voronoi diagrams as convex hull
problems.  Special cases include classical Voronoi diagrams, power
diagrams, order $k$ and farthest point diagrams, Apollonius diagrams,
medial axes, and generalized Voronoi diagrams whose sites are
combinations of points, spheres and half-spaces.  We describe the
application of these results to algorithms for computing 
generalized Voronoi diagrams and find the complexity of
these algorithms.
\end{abstract}


\maketitle

\section{Introduction}\label{sec: introduction}

\subsection{History}
The classical Voronoi diagram for a set of points $\calP = \{p_1, p_2,
\ldots, p_n\}$ in the Euclidean plane is the subdivision of the plane
into Voronoi cells, one for each point in $\calP$.  The Voronoi cell for
the point $p_k$ is the set
\[ \Vor(p_k) = \{ x \in \boldR^2 \, | \, \text{$\dist(x,p_k) \le
\dist(x,p_i)$ for all $i$} \} \] of points $x$ in the plane so that
$p_k$ is a point in $\calP$ closest to $x$.  This notion is so fundamental
that it arises in a multitude of contexts, both in theoretical
mathematics and in the real world.  See \cite{aurenhammer1991voronoi,
deBerg1997, devadoss2011discrete, Aurenhammer2013} for background and
surveys.

The notion of Voronoi diagram may be expanded by changing the
underlying geometry, by allowing the sites to be sets rather than
points, by weighting sites, by subdividing the domain based on
farthest point rather than closest point, or by subdividing the domain
based on which $k$ sites are closest. 
Because there are so many applications of Voronoi diagrams in
applications, it is of considerable interest to find efficient
algorithms for their computation.

By representing Voronoi cells in $\boldR^d$ as projections of the
facets of certain polyhedra in $\boldR^{d+1}$, Brown was able to
use geometric inversion to translate the computation of Voronoi cells
into a convex hull problem that can be implemented with efficient
algorithms~\cite{brown1979voronoi}.
Aurenhammer and Edelsbrunner provided rigorous justification for
Brown's work in \cite{aurenhammer1984optimal}.  Edelsbrunner and
Seidel \cite{edelsbrunner-seidel-86} and others have translated
various types of generalized Voronoi diagrams to {\em lifting
problems} wherein the diagram in $\boldR^d$ is computed by translating
quadratic defining conditions in $\boldR^d$ to linear conditions in
one dimension higher, so that the cells of the diagram are projections
of the facets of a polyhedron in $\boldR^{d+1}$ to $\boldR^d$.
Finding the necessary polyhedron is a {\em convex hull problem}.  In
some cases, generalized Voronoi diagrams are translated to {\em
embedding problems:} the desired diagram in $\boldR^d$ is identified
with a subset of a different generalized Voronoi diagram in a
higher-dimensional Euclidean space. Then the desired
diagram may be found by intersection and projection.  Lifting and
embedding methods can be used to find multiplicatively weighted
Voronoi diagrams \cite{aurenhammer1984optimal}, Voronoi diagrams for
spheres \cite{aurenhammer1987power}, power diagrams
\cite{aurenhammer1987power}, order $k$ Voronoi
diagrams~\cite{rosenberger1991order}, and farthest point Voronoi
diagrams~\cite{brown1979voronoi, aurenhammer2006farthest}.

\subsection{Spaces of spheres}

Viewing Voronoi diagrams and Delaunay triangulations in terms of empty
spheres goes back almost a century \cite{delaunay1934sphere}.  The
classical Voronoi diagram for a set $\calP$ of point sites in $\boldR^d$
may be characterized solely in terms of empty spheres without any
reference to distance functions.  A point $y$ in $\boldR^d$ is in the
Voronoi cell $\Vor(p_k)$ for the point $p_k$ if and only if $y$ is the
center of an empty sphere incident to $p_k$, where ``empty'' means
that the sphere has no point sites from the set $\calP$ in its interior.
The point $y$ is on the boundary of $\Vor(p_k)$ if and only if it is
the center of an empty sphere which has $p_k$ and at least one other
site from $\calP$ incident to it.  Thus, the classical Voronoi cells and
their boundaries can be found from incidence properties of spheres in a space of
empty spheres.  See \cite{devillers1992space} for a mathematical model
of the space of unoriented spheres and applications to various types
of Voronoi diagrams.  

The space of unoriented spheres and hyperplanes lies in the realm of
M\"obius geometry.  We work in the more general setting of Lie sphere
geometry, which is concerned with the moduli space of {\em Lie
spheres}: points (including a point at infinity), oriented spheres,
and oriented hyperplanes.  The basic objects are Lie spheres, so that
oriented hyperplanes, oriented spheres and points are all viewed as
points in the Lie quadric rather than as subsets of $\boldR^d$.  By
considering oriented spheres and hyperplanes, we may discuss oriented
tangency, in addition to incidence, angle, interiority and
exteriority.  The {\em Lie quadric} $\calQ$ in the projective space
$\boldP^{d+2}$ parametrizes the set of Lie spheres.  The Lie quadric
is equipped with a natural symmetric bilinear form (the Lie product)
and a Lie group that acts transitively on the space.  Both Laguerre
geometry and M\"obius geometry are subgeometries of Lie sphere
geometry in the sense of Klein's Erlangen program.

There are several benefits to using Lie sphere geometry for geometric
problems in Euclidean space.  In the Lie quadric, points and spheres are on
the same footing as hyperplanes, and there is no conceptual difference
between inward and outward orientation of spheres.  In this framework,
the set of sites for a generalized Voronoi diagram can be a set of points,
or spheres, or hyperplanes, or a mix of these, with no real
computational difference, thus housing many different kinds of
diagrams under a single roof.  Another benefit is that second order
geometric conditions on spheres in Euclidean space such as incidence,
angle and tangency become first order conditions on points in $\calQ$
in terms of the Lie inner product.  The solution of the geometric
conditions becomes a convex hull problem, yielding a polyhedron.  The
cost of the linearity is that Lie quadric $\calQ$ of $\boldP^{d+2}$ is
defined with a second order condition, so mapping solutions back to
$\boldR^d$ may require intersecting a polyhedron with a quadric.

\subsection{Overview}
We view generalized Voronoi diagrams as minimization diagrams for sets
of functions where each function measures some kind of generalized distance to some 
object. We consider two classes of generalized
Voronoi diagrams in $\boldR^d$. Both are minimization diagrams for
functions that can be expressed in terms of the Lie quadric and the
Lie inner product.

The first class is generalized Voronoi diagrams which can be
formulated in terms of extremal spheres.  The second class is of
minimization functions on $\boldR^d$ which may be expressed in terms
of affine functions with domain $\boldR^{d+3}$ restricted to standard
coordinates for the
Lie quadric.  These two classes include infinitely many generalized
Voronoi diagrams encompassing many classes previously studied, along
with new classes.

In Section \ref{section: main results}, we formally state the main
theorems, Theorem \ref{thm: empty sphere} and Theorem \ref{thm: affine
on quadric}, and in Theorem \ref{thm: complexity} we  show that the
efficiency of an algorithm for computating generalized Voronoi
diagrams with defining data as in the hypotheses of Theorem \ref{thm:
empty sphere} is $O(n\log{n} +
n^{\lfloor\frac{d+3}{2}\rfloor} + nd^2)$.
Theorem \ref{thm: empty sphere} involves generalized
Voronoi diagrams defined in terms of empty spheres, and Theorem 
\ref{thm: affine on quadric} is about minimization diagrams for
functions that arise via restrictions of affine functions to submanifolds.
Preliminaries on minimization diagrams and Lie sphere geometry are in
Section \ref{section: preliminaries}.  Section \ref{section: extremal
spheres} is dedicated to the proof of Theorem \ref{thm: empty sphere}
and applications of that theorem. The proof of the   Theorem \ref{thm:
affine on quadric} and its applications are in Section \ref{section:
affine functions}, along with a useful basic result about
minimization diagrams for 
restrictions of functions to images of embeddings, Proposition \ref{prop:
  MD restriction}.  Finally, in Section \ref{section: algorithm}, we
describe the algorithm for computating generalized Voronoi diagrams
with defining data as in Theorem \ref{thm: empty sphere} and we prove
Theorem \ref{thm: complexity}. 

\section{Main results}\label{section: main results}
We view {\em generalized Voronoi diagrams} in $\boldR^d$ as {\em
minimization diagrams} for collections $\calF =\{f_i\}$
of functions from $\boldR^d$ to $\boldR$, where the functions $f_i$ in $\calF$ measure some kind
of generalized distance to some general notion of site.  The
minimization diagram $\MD(\calF)$ for $\calF$ is the set of {\em minimization
regions}, where a minimization region $MR(I)$ for a set of indices $I$
is equal to the closure of the set of $\bfx \in \boldR^d$ so that
$f_i(\bfx) = f_j(\bfx)$ for all for all $i,j \in I$ and $f_i(\bfx) <
f_k(\bfx)$ for all $i,j \in I$ and $k \not \in I$. See Section
\ref{subsection: minimization diagrams} for precise definitions and
details.

We state the two main theorems and describe some of their applications. 
\subsection{Diagrams defined   by  extremal empty  spheres
}

First we establish notation and terminology.  Let $S(\bfq,r)$ denote
the sphere in $\boldR^d$ with center $\bfq$ and positive radius $r.$
The set $\setext(S(\bfq,r)) = \{ \bfx \in \boldR^d \, : \, \| \bfx -
\bfq \| > r\}$ is the {\em exterior of $S(\bfq,r)$} and
$\setint(S(\bfq,r)) =\{ \bfx \in \boldR^d \, : \, \| \bfx - \bfq \| <
r\}$ is the {\em interior of $S(\bfq,r)$}.  For a unit vector $\bfn$
in $\boldR^d$ and a real number $h,$ the {\em closed half-space
defined by $\bfn$ with the height $h$} is the set $H_{\bfn, h} = \{
\bfx \in \boldR^d\, : \, \bfx \cdot \bfn \ge h\}$.  The {\em M\"obius
scalar product} $\rho$ of spheres $S_1 = S(\bfq_1,r_1)$ and
$S_2=S(\bfq_2,r_2)$ in $\boldR^d$ is
\begin{equation}\label{eqn: Mobius}
  \rho(S_1,S_2) = \smallfrac{1}{2} (r_1^2 + r_2^2 - \| \bfq_1 -
\bfq_2\|^2) .\end{equation}

In Theorem \ref{thm: empty sphere}, we consider generalized Voronoi
diagrams whose minimization regions may be described by geometric conditions on spheres in
$\boldR^d$. Different geometric conditions encode different types of Voronoi
diagrams.  For example, as described previously, for the classical
Voronoi diagram with point sites in a set $\calP$, the cell for a
point site $\bfp$ from $\calP$ can be characterized as the set of all
points $\bfx$ satisfying the geometric condition ``there exists a
sphere $S(\bfx,r)$ about $\bfx$ which is incident to $\bfp$ and whose
interior contains no points from $\calP$.''  The sites of the
generalized Voronoi diagrams for Theorem \ref{thm: empty sphere} will
be encoded in a data set with five types of data points.  These data
points play the role of point sites in the classical case.
The data set is a quintuple $(\calP^-, \calP^+, \calH, \calS, \calT)$
of two sets of points, one set of half spaces, and two sets of
spheres.  For each of the five types of data point, there is a
different type of geometric cell-defining condition involving points being centers
of spheres with certain kinds of properties.  The exact forms of those
conditions are presented in the next definition.

\begin{defn}\label{defn: data set}
Let $\calP^-$ and $\calP^+$ be finite sets of points in $\boldR^d,$
let $\calH$ be a finite set of closed half-spaces in $\boldR^d$, and let
$\calS$ and $\calT$ be finite sets of spheres in $\boldR^d,$ where not
all of these sets are empty.  We call the list $\calD =
(\calP^-,\calP^+, \calH, \calS, \calT)$ a {\em data set}.  We call the
elements of $\calP^-\cup \calP^+ \cup \calH \cup \calS \cup \calT$ the
{\em sites} of the data set $\calD$.

 The data set defines a set $\Sph(\calD)$ of   spheres as follows.
A sphere $S$ is in the set $\Sph(\calD)$ if and only $S$
simultaneously satisfies the following closed conditions:
\begin{enumerate}
\item{$\calP^- \subseteq  \overline{\setext(S)} $,}\label{part: point outside}
\item{$\calP^+\subseteq  \overline{\setint(S)}$,}\label{part: point
    inside}
\item{If $\calH \ne \emptyset$, then $S$ is a subset of the polyhedron
$\cap_{H \in \calH}\overline{H}$,}\label{part: in polyhedron}
\item{$\rho(S,T) \le 0$ for all spheres $T$ in $\calS$, and
}\label{part: power}
\item{If $\calT \ne \emptyset$, then $S$ is a subset of the closed set
$\cap_{T \in \calT} \overline{\setext(T)}$  that is the intersection
of the closures of the exteriors of the spheres $T$ in
$\calT$.}\label{part: outside-sphere}
\end{enumerate}
The data set $\calD$ also defines the subset $Z(\calD)$ of sphere
centers in $\boldR^d$.  A point $\bfx$ in $\boldR^d$ is in $Z(\calD)$
if and only if there is a sphere $S$ whose center $\bfx$ is in
$\Sph(\calD).$
\end{defn}

We will see in Corollary \ref{cor: classical Voronoi} that taking all
the sets in the data set to be trivial except for $\calP^-$ yields the
classical Voronoi diagram, while in Corollary \ref{cor: farthest point}
we see that taking all the sets to be trivial except
for $\calP^+$ yields the farthest point Voronoi diagram.  Corollary
\ref{cor: spherical sites} shows that the  Voronoi diagram for spherical
sites arises from  all the sets in the data set being trivial except for $\calT,$
and in Corollary \ref{cor: power diagram 1}, we see that when $\calS$
is the only nontrivial set in the data set, we obtain the power
diagram.  Corollary \ref{cor: medial axis} shows that the medial axis
of a polyhedron results when $\calH$ is the only nontrivial set.
Thus, our notion of generalized Voronoi diagrams defined by a data set
includes familiar diagrams as special cases while also now allowing
combinations of sites of different types.

Our first main theorem, Theorem \ref{thm: empty sphere}, describes the
structure of the minimization regions for generalized Voronoi diagrams
defined by data sets and conditions on spheres as in Definition
\ref{defn: data set}.  We will see in Subsection \ref{space of
spheres} that oriented spheres in $\boldR^d$ can be parametrized by
the set
\[ \Sph(\boldR^{d}) = \{ \bfx \in \boldR^{d+3} \, : \, x_1 + x_2 =
  1,
   \la \bfx, \bfx \ra = 0
   \},\]
 where $ \la \bfx, \bfx \ra = -x_1^2 + x_2^2 + \cdots + x_{d+2}^2 -
 x_{d+3}^2.$
 This is a subset of the cone
\[Q^d = \{ (x_1, \ldots, x_{d+3}) \, : \,
  -x_1^2 + x_2^2 + \cdots + x_{d+2}^2 - x_{d+3}^2 = 0\} \subseteq
  \boldR^{d+3}.\]
The projection map $\prc: (x_1, \ldots, x_{d+3}) \mapsto (x_3, \ldots,
x_{d+2})$ sends a sphere in $\Sph(\boldR^{d}) $ to its center in
$\boldR^d.$ Now we are ready to state the theorem. 

\begin{thm}\label{thm: empty sphere}
 Let $\calP^-$ and $\calP^+$ be finite sets of points in $\boldR^d,$
 $\calH$ be a finite set of half-spaces in $\boldR^d$, and let
$\calS$ and $\calT$ be finite sets of spheres in $\boldR^d$, where not
all of these sets are empty.  Let $\Sph(\calD)$ be the set of spheres
determined by the data set $\calD=(\calP^-,\calP^+, \calH, \calS,
\calT)$ and let $Z(\calD)$ be the subset of sphere centers in
$\boldR^d$ determined by the data set $\calD$.  Assume that $\Sph
(\calD)$ is nonempty.
Then
\begin{enumerate}
\item Each site in the data set can be assigned a non-strict
homogeneous linear inequality in $\boldR^{d+3}.$ The solution to this
system of inequalities is 
a nonempty polyhedron $P$  in $\boldR^{d+3}$. \label{part: one}
\item   $P \cap \Sph(\boldR^{d}) $ parametrizes the set of spheres
$\Sph(\calD)$ defined by the data set $\calD$.
  \item The boundary of $P \cap \Sph(\boldR^{d}) $ corresponds to the
set of spheres so that at least one condition from Definition
\ref{defn: data set} hold extremally.  The number of conditions which
hold extremally is equal to the number of inequalities from Part 
\eqref{part: one} for which
equalities hold.
  \item The set of sphere centers $Z(\calD)$ is equal to the image
$\prc (P \cap \Sph(\boldR^{d}) )$ of $P \cap \Sph(\boldR^{d})$ under
the projection map $\prc$ from $\boldR^{d+3}$ to $\boldR^d$.
\end{enumerate}
\end{thm}

The inequalities defining the polyhedron from the
data set $\calD = (\calP^-,\calP^+, \calH, \calS, \calT)$ can be found
in 
Table \ref{table: inequalities}.  In Corollaries \ref{cor: classical
Voronoi}, \ref{cor: farthest point}, \ref{cor: spherical sites}, \ref{cor:
power diagram 1}, and \ref{cor: medial axis} to Theorem \ref{thm:
empty sphere}, we obtain the previously known results that the
problems of computing classical Voronoi diagram, the farthest point
Voronoi diagram, the Voronoi diagram for spherical sites, the power diagram, and the
medial axis of a convex polygon are ``convex hull problems''; that is,
the diagrams may be obtained via the intersection of half-spaces in
higher dimensions.  See Table \ref{table: overview}.  An example of a
new type of diagram that can be found with the theorem is in Example
\ref{example: mixed sites}.
\begin{table}
    \setlength{\extrarowheight}{2pt}
   \centering
\begin{tabular}{| c | c | c |}
 \hline
  Diagram type&  \makebox[.15\textwidth]{Inequality} \\
 \hline
 \hline
 Classical Voronoi diagram with data $\{ \xi_i\}$
             & $\la \xi_i, \sigma \ra \le 0 $ \\
 \hline  Voronoi diagram with spherical sites  $\{
  \sigma_i\}$& $\la \sigma_i, \sigma \ra \le 0 $ \\
 \hline
 Farthest point diagram with data $\{ \xi_i\}$&$\la \xi_i, \sigma \ra \ge 0 $  \\
 \hline
 Order two Voronoi diagram with data $\{ \xi_i\}$& $\la \xi_i +  \xi_j, \sigma \ra \le 0 $ \\
 \hline
  Power diagram with data $\{ \sigma_i\}$&
    $ \la \sigma_i', \sigma \ra \le 0 $ \\
 \hline
 Medial skeleton with data $\{ \pi_i\}$& $\la \pi_i, \sigma \ra \le 0 $ \\
 \hline
\end{tabular}
\smallskip
   \caption{Generalized Voronoi diagrams defined by incidence or
     tangency conditions on oriented spheres.  The  notation is
     defined in  Section \ref{section: preliminaries}.}
   \label{table: overview}
 \end{table}

\begin{remark}
  Some conditions on the centers and radii of spheres
  in Theorem \ref{thm: empty sphere}, such as the
center of the sphere lying in a half-space or the radius having a
fixed value or bound, may be expressed in linear equalities or
inequalities in the variables $\sigma_1, \ldots, \sigma_{d+3}$.  In
order to simplify the notation and exposition, we have not
incorporated these options into the statement of Theorem \ref{thm:
empty sphere} and leave this as a simple exercise for the reader.
\end{remark}

\subsection{Diagrams for restrictions of affine functions}
 
 In Section \ref{section: affine functions}, we consider minimization
diagrams for families of functions defined on subsets of $\boldR^d$
that arise from affine functions with domain $\boldR^{d+3}.$ To be
precise, let $\phi: \boldR^d \to \boldR^{d+3}$ be the embedding of
$\boldR^d$ into the previously defined cone $Q^d$  with 
\[  \phi(\bfx) = \left( \smallfrac{1}{2} (1+\bfx \cdot \bfx),
    \smallfrac{1}{2}(1 - \bfx  \cdot   \bfx), \bfx, 0\right)\]
  for $\bfx$ in $\boldR^d.$ Let $L: \boldR^{d+3} \to \boldR$ be an
affine function.  Associate to $L$ the function $f_L = L \circ \phi$
from $\boldR^d$ to $\boldR.$  Theorem \ref{thm: affine on quadric} says
that the minimization diagram  for families  of function of form $f_L$
may be found by solving linear systems of inequalities in
$\boldR^{d+3}$ and pulling the faces of the resulting polyhedron back
to $\boldR^d$ using $\phi.$ 

\begin{thm}\label{thm: affine on quadric}
Let $\calL = \{ L_i\}_{i=1}^n$ be a family of affine functions from
$\boldR^{d+3}$ to $\boldR.$ Let $\calF = \{f_{L_i}\}$ be the quadratic
functions $f_{L_i} = L_i \circ \phi$ from $\boldR^{d}$ to $\boldR$ as defined above.  Let
$\MD(\calF) = \{ \MRF(I) \}$ denote the minimization diagram
for $\calF$, where $\MRF(I) \subseteq \boldR^{d}$ is the minimization region
for indices $I \subseteq [n]$;  and let $\MD(\calL) = \{
\MRL(I) \}$ denote the minimization diagram for the family of
affine functions $\calL$, where $\MRL(I) \subseteq
\boldR^{d+3}$ is the minimization region for the index set $I \subseteq [n]$. 
The minimization diagram for $\calF$ has minimization regions 
\[ \MRF (I) = \phi^{-1}(\MRL (I) \cap Q^d), \quad
  \text{for all 
 $ I \subseteq [n]$}.\]
In particular, each minimization region is the pre-image under $\phi$ of a
facet of a polyhedron and the quadric hypersurface $Q^d.$
\end{thm}

It is natural to ask which functions $f: \boldR^d \to \boldR$ are of
the form $L \circ \phi$ as in the previous theorem.  The next
proposition answers that question.
\begin{prop}\label{prop: of form}
 A function $f: \boldR^d \to \boldR$  is of form $L \circ \phi$ if and
 only if it is of form 
  \[  f(\bfx) = a \| \bfx - \bfq \|^2 + \bfb \cdot \bfx + c, \enskip
    \text{where $a, c \in \boldR$ and $\bfb \in \boldR^{d}$.} \]
\end{prop}
It follows that the following kinds of diagrams in $\boldR^d$ may be
expressed in terms of polyhedra in $\boldR^{d+3}$ as in Theorem
\ref{thm: affine on quadric}: classical Voronoi diagrams, farthest
point Voronoi diagrams, power diagrams, additively weighted Voronoi
diagrams, multiplicatively weighted
Voronoi diagrams, and medial axes of convex polygons.
(See Proposition \ref{prop: affine diagrams}.)  

\subsection{Order $k$ generalized Voronoi diagrams}
Let $\calF = \{f_1, f_2, \ldots, f_m\}$ be a set of continuous
functions mapping $\boldR^d$ to $\boldR.$  Let $k \in [m]$.  The {\em
order $k$ Voronoi diagram $\Vor^k(\calF)$ for $\calF$} is the set
of cells indexed by $k$-combinations of $[m]$, where
each cell is determined by the $k$
smallest values among the values of $f_i(\bfx)$, with $i \in [m].$
The {\em cell} $\Cell(i)$ for the subset $i = \{i_1, i_2,
\ldots, i_k\}$ of size $k$ is the closure of the set of all $x$ in
$\boldR^d$ so that
\[ f_{j}(\bfx) \le f_{l}(\bfx)  \quad \text{for $j \in \{i_1, i_2, \ldots, i_k\}$ and
       $l \not \in \{i_1, i_2, \ldots, i_k\}$}.\]
The cells for the order 1  Voronoi diagram $\Vor^1(\calF)$ are a
subset of  the usual minimization diagram $\MD(\calF)$ for $\calF.$

\begin{prop}\label{prop: order k}
  Let $\calF = \{f_i \}$ be a family of continuous real-valued functions with
domain $\boldR^d$, where for each $i,$
  \[ f_i(\bfx) = a_i \| \bfx - \bfq_i \|^2 + \bfb_ i\cdot \bfx + c_i,
\enskip \text{where $a_i, c_i \in \boldR$ and $\bfb_i \in
  \boldR^{d}$.} \]
The order $k$ Voronoi diagram $\Vor^k(\calF)$  for $\calF$ may be
computed as in Theorem \ref{thm: affine on quadric}.  
\end{prop}

Specifically, the cells of $\Vor^k(\calF)$ are minimization regions
for a minimization diagram defined by a family of functions
$\calF^k$.  See Proposition \ref{prop: order k bijection}.
This is  similar to Theorem 4.8 of \cite{boissonnat2018geometric}, which describes the
classical order $k$ Voronoi diagram for a set $S$ of points in
$\boldR^d$  as a weighted Voronoi diagram for
which the points are the centroids of cardinality $k$ subsets of $S$.

\subsection{Computation}
Suppose that we have a data set $\calD=(\calP-, \calP^+, \calH, \calS,
\calT)$ of points, hyperplanes, and hyperspheres as in Definition
\ref{defn: data set} defining a polytope as in Theorem \ref{thm: empty
sphere}.  Let $\MD (\calD)$ be a minimization diagram defined by the
data set.  we compute the minimization diagram by finding the
halfpsace intersection and intersection with the Lie quadric, then
mapping back to $\boldR^d$ to get the minimization diagram.  The
algorithm is given in Section \ref{section: algorithm}.

\begin{thm}\label{thm: complexity} The algorithm in Section
  \ref{section: algorithm} is $O(n\log{n} +
n^{\lfloor\frac{d+3}{2}\rfloor} + nd^2)$.
\end{thm}

\begin{acknowledgments}
  We are grateful to Idaho State University's Student Career Path
Internship program and Office of Research for funding Egan Schafer's
work on this project while he was a student.
This project was completed during the second
author's sabbatical visit at the University of Arizona.  She is
grateful to Dave Glickenstein for hosting her visit and for 
many illuminating discussions on various aspects of computational geometry.
\end{acknowledgments}

\section{Preliminaries}\label{section: preliminaries}

In this section we review minimization diagrams and Lie sphere
geometry.  
\begin{conventions} We will always work in $\boldR^d,$ with $d \ge 2$.
  We denote the standard orthonormal basis vectors in $\boldR^d$ by
$\bfe_1, \bfe_2, \ldots, \bfe_d$, where the entries of $\bfe_k$ are
all zero, except the $k$th entry, which is a one.  We represent the
coordinates of a vector $\bfx = (x_i)$ in $\boldR^d$ by $x_1, \ldots,
x_d.$ For $n \in \boldN,$ we denote the set $\{1, 2, \ldots, n\}$ by
$[n].$

When we work with $d$-dimensional projective space $\boldP^d$, we
always identify $\boldP^d$ with nonzero points in $\boldR^{d+1}$
modulo rescaling: $\boldP^{d} = \{ [\bfx] \, : \, \bfx \in
\boldR^{d+1} \setminus \{ \bfzero \} \}$, where $[\bfx]$ is the
equivalence class $\{ \lambda \bfx \, : \, \lambda \ne 0 \}$ for
any nonzero $\bfx$ in $\boldR^{d+1}$.

We always assume that the data defining a generalized Voronoi
diagram or minimization diagram consists of two or more sites or
functions.
\end{conventions}

\subsection{Minimization diagrams}\label{subsection: minimization diagrams}
We view generalized Voronoi diagrams in $\boldR^d$ as minimization
diagrams.   Let $\calF = \{f_1, f_2, \ldots, f_m\}$ be a set of continuous
real-valued 
functions with domain $\boldR^d$.  The {\em lower envelope of
$\calF$} is the function $f^-: \boldR^d \to \boldR$ defined by $
f^-(\bfx) = \min \{ f_i(\bfx) \, : \, i \in [m] \}.$ For $\bfx\in \boldR^d,$
the {\em minimal index set $I(\bfx)$ \em for $\bfx$} is the set of indices $i$ for
which the function values $f_i(\bfx)$ of $f$ at $\bfx$ are minimal:
\begin{equation}\label{defn I}
  I(\bfx) =\{ i \in [m] \, : \, f_i(\bfx) = f^-(\bfx) \}.\end{equation}
For a subset $J$ of $[m],$ the {\em minimization region $MR(J)$ for
$J$} is the set of the $\bfx$ in $\boldR^d$ so that the minimum value of
$\{ f_i(\bfx) \, : \, i \in [m] \}$ is achieved by precisely the
functions with indices in $ J$:
\begin{equation}\label{defn MR}
  MR(J) =  \{ \bfx \in \boldR^d \, : \, I(\bfx) = J\}. \end{equation}
For $i \in [m]$, the {\em minimization cell} $\Cell(i)$ is defined by 
\[ \Cell(i) = 
\{ \bfx \in \boldR^d \, : \, f_i(\bfx) \le f_j(\bfx) \enskip
\text{for  all $j \in [m]$}         
\}. \] 
It is the union of all minimization regions $MR(J)$ with $i\in J.$
The {\em minimization diagram $ \MD(\calF)$ for $\calF$} is the
collection $\{ \overline{MR(J)} \, : \, J \subseteq [m]\}$ of the
closures of all minimization regions.

\begin{example}\label{example: classical Voronoi min diagram}
  Given a set $\calP= \{p_1, \ldots, p_m \}$ of point sites in
$\boldR^d,$ if $f_i(\cdot) = \dist(\cdot, p_i)$ for $i \in [m]$,
then the minimization diagram for $\calF= \{f_i\}_{i=1}^m$ yields the
classical Voronoi diagram for the set $\calP$.  The lower envelope
$f^-$ measures the distance from a point in in $\boldR^d$
to the set $\calP$ and the minimal index set $I(\bfx)$ for $\bfx$
in $\boldR^d$ consists of the indices of the point or points in $\calP$
that are closest to $\bfx.$ For $J \subseteq [m]$, the minimization
region $MR(J)$ is the set of points $\bfx$ whose set of closest points
is $\{p_j \, : \, j \in J\}.$ The classical Voronoi cells are the
closures of the minimization cells
$\Cell(\{i\}) = \{ \bfx \in \boldR^d \, : \, \dist(\bfx,p_i) \le \dist (\bfx,
p_j) \enskip \text{for all $j \in [m]$} \} $.  
\end{example}

  The {\em power} $\Pow(\bfp,S(\bfq,r))$ of a point $\bfp$ in
$\boldR^d$ with respect to a sphere $S(\bfq,r) \subseteq \boldR^d$
with center $\bfq$ and radius $r$ is
\begin{equation}\label{eqn: power def}
  \Pow(\bfp,S(\bfq,r)) = \dist^2(\bfp,\bfq) - r^2.\end{equation}

\begin{example}\label{example: power diagram}
Let $\calS = \{S_1, S_2, \ldots, S_m\}$ be a set of spheres in
$\boldR^d$.  The {\em power diagram of $\calS$} is the minimization
diagram for the set of power functions $\{ \Pow( \cdot,
S_i)\}_{i=1}^m$.
\end{example}

\subsection{M\"obius geometry, Laguerre geometry and Lie sphere geometry}

In this section, we summarize the basics of Lie sphere geometry,
following the reference \cite{cecil08}, omitting the motivation and  proofs that may be found there.

\subsubsection{M\"obius geometry}
For a point $\bfq$ in $\boldR^d$ and a real number $r$ which is
possibly nonpositive, let $S(\bfq,r)$ denote the {\em sphere with center
$\bfq$ and radius $r$}
\[ S(\bfq,r) = \{ \bfx \in \boldR^d \, : \, \| \bfx - \bfq \|^2 =
  r^2\}.\]
If the sphere is oriented, it is endowed with unit normal vector
$\bfn(\bfx) = \frac{1}{r} (\bfq - \bfx)$ when $r$ is nonzero.  When
$r$ is positive, this normal vector is inward, when $r$ is negative,
the normal vector is outward. When $r$ is zero, the sphere is 
an unoriented point sphere.

For unit vector $\bfn$ in $\boldR^d$ and a real number $h,$ the
{\em oriented hyperplane $H_{\bfn, h}$ with normal vector $\bfn$ and
height $h$} is the set
\[ H_{\bfn, h} = \{ \bfx \in \boldR^d \, : \,  \bfx \cdot \bfn = h
  \}. \]
The {\em positive (open) half-space defined by $ H_{\bfn, h} $} and the {\em
negative (open) half-space $H^-_{\bfn, h} $ defined by $ H_{\bfn, h} $} are
\begin{align*}
  H^+_{\bfn, h}& = \{\bfx \in \boldR^d \, : \, \bfx \cdot \bfn >
  h\} \quad \text{and} \\
H^-_{\bfn, h}&=
\{\bfx \in \boldR^d \, : \, \bfx \cdot \bfn < h\},\end{align*}
respectively.  The normal vector points in the direction of increasing
$h$ and ``into'' $H^+_{\bfn, h}.$

Stereographic projection $\rho: \boldR^d \to \boldR^{d+1}$ maps
$\boldR^d$ homeomorphically and conformally onto the punctured unit
sphere $S(\bfzero,1) \setminus \{-\bfe_1\}$ in $\boldR^{d+1}$.  For a
point $\bfx$ in $\boldR^d,$ its image in $ \boldR^{d+1}$ under
stereographic projection is
\[ \rho(\bfx) = \left( \frac{1 - \| \bfx \|^2}{1 + \| \bfx \|^2} ,
    \frac{2}{1 + \| \bfx \|^2} \, \bfx \right) \in \boldR^{d+1}.  \]
The south pole $-\bfe_1$ in $S(\bfzero,1)$ corresponds
to the point at infinity in the one point compactification of
$\boldR^d.$

Let $\boldP^{d+1}$ denote $(d+1)$-dimensional projective space.  Let
$\phi: \boldR^{d+1} \to \boldP^{d+1}$ be the affine embedding of
$\boldR^{d+1}$ into $\boldP^{d+1}$ given by $\phi(\bfx) = [(1,\bfx)]$
for $\bfx \in \boldR^{d+1}$.  The composition $\phi \circ \rho:
\boldR^d \to \boldP^{d+1}$ maps $\bfx$ in $\boldR^d$ to $(\phi \circ
\rho) (\bfx) = [\xi_{\bfx}']$ in  $\boldP^{d+1},$  where 
\begin{equation}\label{eqn: Lie point}
 \xi_{\bfx}' =   \left( \frac{1 + \| \bfx \|^2}{2},  \frac{ 1 - \| \bfx \|^2}{2},
     \bfx \right).\end{equation}
 The range of $\phi \circ \rho$ is
 \[ (\phi \circ \rho) (\boldR^d) = \{ [(1, \bfx)] \,:\, \| \bfx \| = 1
   \} \setminus \{ [(1,-1, \bfzero)]\} \subseteq \boldP^{d+1}.\]
The missing point $[(1,-1, \bfzero)]$ corresponds to the point at
infinity in $\boldR^d$ and is called the {\em improper point}.

Endow $\boldR^{d+2}$ with the Lorentz inner product
\[ B(\bfx,\bfy) = -x_1y_1 + x_2y_2 + \cdots + x_{d+2} y_{d+2}. \]
A point $\bfx$ in $\boldR^{d+2}$ is called {\em timelike}, {\em
spacelike} or {\em lightlike} if $B(\bfx,\bfx)$ is negative, positive
or zero respectively.  Although $B$ is not well-defined on
$\boldP^{d+1}$, it still makes sense to talk about the vanishing,
positivity and negativity of $B([\bfx],[\bfx])$, and the vanishing or
nonvanishing of $B([\bfx],[\bfy])$, for $[\bfx]$ and $[\bfy]$ in
$\boldP^{d+1}$, as these properties are invariant under rescaling by
nonzero scalars.  We say that $[\bfx] \in \boldP^{d+1}$ is {\em
timelike}, {\em spacelike} and {\em lightlike} if $B(\bfx,\bfx)$ is
negative, positive or zero respectively.  Given $[\bfx] \in
\boldP^{d+1},$ let $ [\bfx]^\perp = \{ [\bfy] \in \boldP^{d+1} \, : \,
B(\bfx,\bfy) = 0\}$.

The {\em M\"obius sphere} $\Sigma^d \subseteq \boldP^{d+1}$ is the
$d$-dimensional set
\[
  \Sigma^d = \{ [\bfx] \in \boldP^{d+1} \, : \, B(\bfx,\bfx) = 0\} \]
of lightlike points in $\boldP^{d+1}$.  Since
\begin{align*}
  \Sigma^d  
  &= \{ [\bfx] \in \boldP^{d+1} \, : \, x_1^2 = x_2^2 + \cdots +
    x_{d+2}^2  \} \\
  &= \phi(S(\bfzero,1)) \\
&= (\phi \circ \rho) (\boldR^d)  \cup \{[(1,-1, \bfzero)]\},\end{align*}
the map $\phi \circ \rho$ is a bijection between $\boldR^d \cup \{
\infty \}$ and the M\"obius sphere.

It can be shown that the set of spacelike points in $\boldP^{d+1}$
parametrizes the set of hyperspheres and hyperplanes in $\boldR^d$.
A spacelike point $[\xi ]$ defines the codimension one
subset of the M\"obius sphere:
\begin{equation}\label{eqn: polar} [\xi]^\perp \cap \Sigma^d =
  \{ \zeta \in \Sigma^d \subseteq
  \boldP^{d+1} \, : \, B(\xi,\zeta) = 0\} .\end{equation}
The pre-image of $ [\xi]^\perp \cap \Sigma^d$ under $\phi \circ \rho$
is a hyperplane or hypersphere in $\boldR^d$; which of these it is
depends on whether or not $ [\xi]^\perp \cap \Sigma^d $ contains the
improper point $[(1,-1, \bfzero)]$, or equivalently, whether or not
$\xi_1 + \xi_2$ is zero.  To be precise, the hypersphere $S(\bfp,r)$
is represented by $[\sigma_{\bfp,r}']^\perp\cap \Sigma^d\subseteq
\boldP^{d+1}$, where
\begin{equation}\label{eqn: Lie sphere}
  \sigma_{\bfp,r}' = \left(
    \frac{1 +\bfp \cdot
        \bfp -  r^2}{2},\frac{1 - \bfp  \cdot \bfp + r^2}{2}, \bfp
    \right),\end{equation}
and the hyperplane $H_{\bfn, h}$ is represented by
$[\pi_{\bfn, h}' ]^\perp\cap \Sigma^d \subseteq \boldP^{d+1}$, where
\begin{equation}\label{eqn: Lie hyperplane}\pi_{\bfn,
    h}' = (h, -h, \bfn).\end{equation}
Note that when $r=0$ in Equation \eqref{eqn: Lie sphere}, we get the
the equation for the center of the sphere in Equation \eqref{eqn: Lie
  point}; that is, 
$\sigma_{\bfp,0}^\prime= \xi_{\bfp}'.$ 

The vectors $\xi_{\bfx}', \sigma_{\bfp,r}'$ and $\pi_{\bfn,h}'$ in
$\boldR^{d+2}$ as in Equations \eqref{eqn: Lie point}, \eqref{eqn: Lie sphere} and \eqref{eqn: Lie hyperplane}
are {\em standard coordinates} for the
equivalence classes $[\xi_{\bfx}'], [ \sigma_{\bfp,r}']$ and
$[\pi_{\bfn,h}']$ denoting points in projective space, and we say the
points are in {\em standard form}.  See Table
\ref{table: standard coordinates} for a summary of the one-to-one
correspondence 
between points, spheres and hyperplanes in $\boldR^d \cup \{ \infty
\}$ and subsets of the M\"obius sphere $\Sigma^d$ expressed in standard coordinates.
\begin{table}
    \setlength{\extrarowheight}{2pt}
\begin{center}
\begin{tabular}{| l | c|}
                    \hline
 Object in $\boldR^d \cup \{\infty\}$&
                                                        Counterpart in
                                       $\Sigma^d \subseteq \boldP^{d+1}$   \\
  \hline
  \hline
  Point  $ \bfx$ & {\tiny $\D{  \left[\left(\frac{1 + \bfx \cdot
                 \bfx}{2}, \frac{1 - \bfx \cdot \bfx}{2}, \bfx
                 \right)\right]   }$}  \\
  Improper point  & $\D{ [(1,-1,\bfzero)] }$  \\
  Oriented hypersphere $S(\bfp,r)$ & 
                                                   {\tiny   $\D{\left[\left(\frac{1 +
                                             \bfp \cdot \bfp -
                                             r^2}{2},\frac{1 - \bfp
                                             \cdot \bfp + r^2}{2},
                                             \bfp\right)\right]^\perp}$}$\cap \Sigma^d$\\
  Oriented hyperplane $H_{\bfn, h}$ & $\D{[(h, -h, \bfn)]^\perp \cap \Sigma^d}$  \\
                              \hline
\end{tabular}
\end{center}
\smallskip
\caption{Points,  
  unoriented spheres and unoriented
  hyperplanes  in $\boldR^d \cup \{ \infty \}$, viewed in the M\"obius sphere $\Sigma^d$, expressed using
  standard  coordinates.}
\label{table: standard coordinates}
  \end{table}

See Table \ref{table: Lorentz} for a summary of the values of the
Lorentz inner product of vectors in standard coordinates representing
points, spheres and hyperplanes.  From these values,
it is not hard
to show that for $\bfx \in \boldR^d$,
\begin{itemize}
\item $\bfx$ is incident to $S(\bfp,r) \subseteq \boldR^d$ if and only
if $[\xi_{\bfx}] \in [\sigma_{\bfp,r}']^\perp$,
\item  $\bfx$ is incident to they hyperplane $H_{\bfn, h}$ if and only
if $[\xi_{\bfx}] \in [\pi_{\bfn,h}]^\perp$,
  \item $\bfx$ is inside the sphere $S(\bfq,r)$ if and only if
$B(\xi_\bfx,\sigma_{\bfq,r})>0$, and
  \item $\bfx$ is inside the half space $H_{\bfn, h}^+$ if and only if
$B(\xi_\bfx, \pi_{\bfn,h})>0$.
\end{itemize}

\begin{table}
    \setlength{\extrarowheight}{2pt}
  \centering
\begin{tabular}{| l | c|}
\hline
{Pair of objects in  $\Sigma^d$} & {Lorentz inner product } \\
\hline
  \hline
 Points $\xi_{\bfp}'$ and $\xi_{\bfq}'$ &  $  - \|  \bfp - \bfq \|^2$\\
  \hline
 Sphere $\sigma_{\bfq, r}'$ and   point $\xi_{\bfp}'$ &  $ r^2 - \|  \bfp - \bfq \|^2$\\
  \hline
 Spheres $\sigma_{\bfq_1, r_1}'$ and  $\sigma_{\bfq_2, r_2}'$
                  &  $ \smallfrac{1}{2} (r_1^2 + r_2^2 - \| \bfq_1 -
                    \bfq_2 \|^2)$ \\
  \hline
 Sphere $\sigma_{\bfq, r}'$ and   hyperplane $\pi_{\bfn,
  h}'$ &  $ \bfq \cdot \bfn - h $\\
  \hline
 Point $\xi_{\bfp}'$ and   hyperplane $\pi_{\bfn,
  h}'$ &  $ \bfp \cdot \bfn - h $\\
   \hline
 Hyperplanes $\pi_{\bfn_1, h_1}'$ and  $\pi_{\bfn_2,
  h_2}'$&  $ \bfn_1 \cdot \bfn_2$\\
  \hline
\end{tabular}
\smallskip
\caption{The values  of the
  Lorentz inner product $B(\tau,\zeta)$ of vectors $\tau$ and
  $\zeta$  in $\boldR^{d+2}$, in standard coordinates, that 
  encode spheres,
  hyperplanes, and points. }
\label{table: Lorentz}
\end{table}

 \subsubsection{Lie sphere geometry}
Thus far, we have been working in the setting of M\"obius geometry,
which is the study of unoriented spheres, angles, and conformal
diffeomorphisms of the M\"obius sphere.  Laguerre geometry is
concerned with oriented spheres and oriented hyperplanes, and the
oriented contact of these.  We now widen our scope and move to the
setting of Lie sphere geometry, which has both M\"obius geometry and
Laguerre geometry as subgeometries.  A {\em Lie sphere} is an oriented
sphere, a point in $\boldR^d \cup \{\infty\}$, or an oriented
hyperplane. These are the fundamental objects in Lie sphere geometry.

Endow $\boldR^{d+3}$ with the {\em Lie product} defined by 
\[ \la \bfx,\bfy\ra = -x_1y_1 + x_2y_2 + \cdots + x_{d+2} y_{d+2} -
  x_{d+3} y_{d+3}, \]
for $\bfx$ and $\bfy$ in $\boldR^{d+3}$.  The {\em Lie quadric}
$\calQ^{d+1}$ is the quadric hypersurface in $(d+2)$-dimensional
projective space $\boldP^{d+2}$ 
    \begin{align*}
      \calQ^{d+1} & = \{ [\bfx] \in \boldP^{d+2}\, : \, \la \bfx, \bfx \ra = 0 \} \\
      &= \{ [(x_1, x_2, \ldots, x_{d+3})] \, : \, -x_1^2 + x_2^2 +
        \cdots + x_{d+2}^2 - x_{d+3}^2 = 0 \}.
    \end{align*}
  
Now we describe a bijection between Lie spheres in $\boldR^d$ and
points in the Lie quadric $\calQ^{d+1}.$ The points in $\calQ^{d+1}$
are called {\em Lie coordinates} for the corresponding Lie spheres as
subsets of $\boldR^d$.  The point $\bfp$ in $\boldR^d$ maps to
$[\xi_{\bfp} ]$ in $\calQ^{d+1}$, where
\begin{equation}\label{eqn: lie point} \xi_{\bfp} = \left(\frac{1 +\bfp \cdot \bfp }{2},\frac{1 - \bfp
      \cdot \bfp}{2}, \bfp, 0\right),\end{equation}
and the point at infinity maps to $[(1,-1,\bfzero, 0)]$.  The oriented
hypersphere $S(\bfp,r)$ maps to $[\sigma_{\bfp,r}]$ in $\calQ^{d+1}$,
where
\begin{equation}\label{eqn: lie sphere}
  \sigma_{\bfp,r} = \left(
    \frac{1 +\bfp \cdot
        \bfp -  r^2}{2},\frac{1 - \bfp  \cdot \bfp + r^2}{2}, \bfp,
      r
    \right),\end{equation}   
  and the oriented hyperplane $H_{\bfn, h}$ maps to $[\pi_{\bfn,h}]$
in $\calQ^{d+1}$, where
\[ \pi_{\bfn,
    h} = (h, -h, \bfn, 1).\]
Recall that a positive radius for a sphere indicates an
inward-pointing normal vector.

 The M\"obius sphere $\Sigma^d$ is embedded in the Lie quadric
$\calQ^{d+1}$ as the set of points with last coordinate equal to zero:
 \[ \Sigma^d \cong [\bfe_{d+3}]^\perp \cap \calQ^{d+1} = \{ [\bfx] \in
\calQ^{d+1} \, : \, \la \bfx, \bfe_{d+3}\ra = 0 \} .  \] 
Furthermore, if a point in this set is not the point at infinity, it
is of form $\left(\frac{1 + \bfx \cdot \bfx}{2}, \frac{1 - \bfx \cdot
\bfx}{2}, \bfx, 0\right)$ and therefore lies in the affine hyperplane $x_1 + x_2 =
1$. (See Table \ref{table: Lie coordinates properties}.)  From the
equation of the quadric we get $-x_1 + x_2 + x_3^2 + \cdots +
x_{d+2}^2= 0$, and if we let $x_0 = x_1-x_2,$ we have $x_0 = x_3^2 +
\cdots + x_{d+2}^2,$ a paraboloid $\calR$ in the intersection of the
hyperplanes $x_{d+3}=0$ and $x_1 + x_2 = 1$.  Thus, the set 
\begin{align}\label{eqn: paraboloid}
  \calR
        &= \left\{
          \left( \frac{1+\bfx \cdot \bfx}{2}, \frac{1 - \bfx  \cdot
                                              \bfx}{2}, \bfx,
                                              0\right) \,: \, \bfx \in
          \boldR^d           \right\} \\ \notag
        &= \left\{ \bfx\in \boldR^{d+3} \, :
          \, x_1 + x_2 = 1,   x_2 - x_1 = x_3^2 + \cdots +
x_{d+2}^2, x_{d+3} = 0  \right\} 
  \end{align}
  parametrizes the set of point spheres in $\Sph(\boldR^d)$.

Note that if we project to the first $d+2$ coordinates, the point
$\xi_{\bfp}$ maps to the point $\xi_{\bfp}'$ as in Equation
\eqref{eqn: Lie point}, the oriented sphere $\sigma_{\bfp,r}$ maps to
the unoriented sphere $\sigma_{\bfp,r}'$ as in Equation \eqref{eqn:
Lie sphere}, and the oriented hyperplane $\pi_{\bfn, h} $ maps to the
unoriented hyperplane $\pi_{\bfn, h}'$ as in Equation \eqref{eqn: Lie
  hyperplane}.  As in  the M\"obius setting, points may be viewed as
spheres of zero radius, as $\sigma_{\bfp,0}=\xi_{\bfp}.$

  When the elements $[\xi_{\bfp}], [\sigma_{\bfp,r}]$ and
$[\pi_{\bfn,h}]$ of $\boldP^{d+2}$ are represented with the points
$\xi_{\bfp}, \sigma_{\bfp,r}$ and $\pi_{\bfn,h}$ in $\boldR^{d+3},$
respectively, the coordinates are called {\em standard Lie
  coordinates}, and we say that the points are in {\em standard form}.
If an oriented sphere $[\bfx] \in \calQ^{d+1}$ is in
standard Lie coordinates, the last coordinate is the {\em signed radius}.
Every element $[\zeta]$ of $\calQ^{d+1}$ can be uniquely represented
in one of these three standard forms, in which $\zeta_1 + \zeta_2 = 1$
as in the first two, or $\zeta_1 + \zeta_2 =0$ and $\zeta_{d+3}=1$ as
in the last.  The bijection between Lie spheres in $\boldR^d$ and
points $[\bfx]$ in $\calQ^{d+1}$ is summarized in Table \ref{table:
Lie coordinates}.  The third and fourth columns in Table \ref{table:
Lie coordinates properties} show how to determine the type of a Lie
sphere based on the values of $x_1+x_2$ and $x_{d+3}$ when it is
expressed in standard coordinates $[\bfx].$ The table shows how to
invert the map given in Table \ref{table: Lie coordinates} and send a
point in the Lie quadric to the corresponding set in $\boldR^d \cup \{
\infty\}$. Thus, the Lie quadric can be seen as the
moduli space of Lie spheres.

A Lie sphere $[\bfx]$ in the Lie quadric corresponds to a point in
$\Sigma^d \sim \boldR^d \cup \{ \infty \}$ if and only if $x_{n+3}= 0$.  If a Lie
sphere $[\bfx]$ has $x_{n+3}\ne 0$, then $\bfx = (\bfw, \pm
\sqrt{B(\bfw,\bfw)}),$ where $[\bfw]$ is a spacelike point in
$\boldP^{d+1}$; the map $[\bfx] \mapsto [\bfw]$ is a double covering
of the spacelike points (hyperspheres and hyperplane) corresponding to
adding orientations.

\begin{table}
    \setlength{\extrarowheight}{4pt}
\begin{center}
\begin{tabular}{| l | c|}
                    \hline
 {Objects in $\boldR^d$} & {Points in $\calQ^{d+1}$ }   \\
  \hline
  \hline
  Point $\bfx$ & {\tiny $\D{ \left[\left(\frac{1 + \bfx \cdot \bfx}{2}, \frac{1 - \bfx \cdot \bfx}{2}, \bfx, 0\right)\right]}$}  \\
  Point at infinity & $\D{[(1,-1,\bfzero, 0)]}$  \\
  Oriented hypersphere $\sigma_{\bfp, r}$  &
                                                   {\tiny   $\D{\left[\left(\frac{1 +
                                             \bfp \cdot \bfp -
                                             r^2}{2},\frac{1 - \bfp
                                             \cdot \bfp + r^2}{2},
                                             \bfp, r\right)\right]}$}\\
  Oriented hyperplane $\pi_{\bfn, h}$ & $\D{[(h, -h, \bfn, 1)]}$  \\
                              \hline
\end{tabular}
\end{center}
\smallskip
    \caption{Points in the Lie quadric $\calQ^{d+1}
      \subseteq \boldP^{d+2}$, in
      standard coordinates, for Lie spheres in $\boldR^d$.}\label{table: Lie coordinates}
  \end{table}
  
\begin{table}
  \begin{center}
      \setlength{\extrarowheight}{2pt}
\begin{tabular}{| l |c|c|c|c|}
                    \hline
  {Objects in $\boldR^d$} &  $x_0 = x_2 - x_1$ & $x_1 + x_2$
  &  $x_{d+3}$   & $(x_3, \ldots, x_{d+2})$\\
  \hline
  \hline
  Point $\bfx$ & $ -\| \bfx \|^2$ &  $1$ & $0$ & $\bfx$ \\
  Point at infinity & $2$ & $0$ & $0$ & $\bfzero$   \\
  Oriented hypersphere $\sigma_{\bfp, r}$  & $-\| \bfp \|^2 +  r^2$& $1$ & $r \ne
                                                                0$ &
                                                                     $\bfp$ \\
  Oriented hyperplane $\pi_{\bfn, h}$ & $2h$ & $0$  & $1$ & $\bfn$ \\
                              \hline
\end{tabular}
\end{center}
\smallskip
\caption{Properties of standard Lie coordinates for Lie
  spheres. }
\label{table: Lie coordinates properties}
\end{table}    
\smallskip

Table \ref{table: Lie product} shows the values of the Lie inner
product for various combinations of Lie spheres expressed in standard
Lie coordinates.   (Note that the distance from a point to a sphere,
$\| \bfp - \bfq \| - r$, does not arise via the Lie inner product.)
Geometric consequences of the values in the table are summarized in
Table \ref{table: Lie properties}.  In particular, 
oriented contact of Lie spheres is now discernable.  Two oriented
hyperspheres, or an oriented hypersphere and an oriented hyperplane
are said to be in {\em oriented contact} if they are tangent at a
point and have the same normal vector at that point.
The Lie inner product of the standard Lie coordinates
$\sigma_{\bfq_1,r_1}$ and $\sigma_{\bfq_2,r_2}$ for oriented spheres
$S(\bfq_1,r_1) \subseteq \boldR^d$ and $S(\bfp_1,r_1) \subseteq
\boldR^d$ is
\begin{equation}\label{equation: lie spheres ip} \la \sigma_{\bfq_1,r_1}, \sigma_{\bfq_2,r_2} \ra
= \smallfrac{1}{2} \left( ( r_1 - r_2)^2 - \| \bfq_1 -
  \bfq_2 \|^2\right) .\end{equation}
This is zero precisely when the spheres are in oriented contact.
Two unoriented spheres $S(\bfq_1,r_1)$ and   $S(\bfq_2,r_2)$ with
positive radii  are tangent with outward normal vectors pointing in
opposite directions if and only if the oriented spheres
$S(\bfq_1,r_1)$ and   $S(\bfq_2,-r_2)$ are in oriented contact.  In
this situation we say that  $S(\bfq_1,r_1)$ and   $S(\bfq_2,r_2)$ are
{\em externally tangent}.

\begin{table}
    \setlength{\extrarowheight}{2pt}
  \centering
\begin{tabular}{| p{3in} | c|}
\hline
{Type of pairing} & {Lie product $\la \cdot , \cdot \ra$} \\
\hline
  \hline
 Points $\xi_{\bfp}$ and $\xi_{\bfq}$ &  $  -  \frac{1}{2} \|  \bfp - \bfq \|^2$\\
  \hline
 Oriented sphere $\sigma_{\bfq, r}$ and   point $\xi_{\bfp}$ &  $
                                                               \frac{1}{2}
                                                               (r^2 - \|
                                                    \bfp - \bfq \|^2)$\\
  \hline
 Oriented spheres $\sigma_{\bfq_1, r_1}$ and  $\sigma_{\bfq_2, r_2}$
                  &  $ \frac{1}{2} \left(( r_1 - r_2)^2 - \| \bfq_1 -
                    \bfq_2 \|^2\right)$ \\
  \hline
 Oriented sphere $\sigma_{\bfq, r}$ and   hyperplane $\pi_{\bfn,
  h}$ &  $\bfq \cdot \bfn - h - r$\\
  \hline Oriented hyperplanes $\pi_{\bfn_1,
  h_1}$ and  $\pi_{\bfn_2,
  h_2}$ &  $\bfn_1 \cdot \bfn_2 - 1$\\
  \hline
\end{tabular}
\smallskip
\caption{The values  of  the  Lie inner
  product of points in the  Lie
  quadric    represented in standard form.  Values  for
  points $\bfp$  may be found using the equality $\xi_{\bfp}=\sigma_{\bfp,0}$.}
\label{table: Lie product}
\end{table}
 The Lie inner product defines a map from $\boldR^{d+3}$ to its dual
space $(\boldR^{d+3})^\ast$ by sending the vector $\bfa$ to the linear
functional $\la \bfa, \cdot \ra.$ This yields an isomorphism between
$1$-dimensional subspaces in $\boldR^{d+3}$ and ($d+2$)-dimensional
subspaces in $\boldR^{d+3}$, and from there, a duality between points
in $\boldP^{d+2}$ and hyperplanes in $\boldP^{d+2}$.  For any points
$[\zeta]$ and $[\tau]$ in projective space $\boldP^{d+2}$ with
corresponding projective hyperplanes $[\zeta]$ and $[\tau]$ in
$\boldP^{d+2}$,
\[  [\zeta] \in [\tau]^\perp \, \Longleftrightarrow \, \la \zeta,
  \tau \ra
  = 0 \,  \Longleftrightarrow \,   [\tau] \in [\zeta]^\perp. \]
Thus, the conditions involving incidence, orthogonality and
oriented contact in Table \ref{table: Lie properties} may now be
interpreted as dual conditions between points and hyperplanes in
projective space.

\begin{table}
  \centering
  \setlength{\extrarowheight}{2pt}
\begin{tabular}{|c | p{2.75in} | c |}
\hline
&  \makebox[2.75in][c]{Property of Lie spheres} &
                                                        Property of corresponding \\
 & \makebox[2.75in][c]{as objects in $\boldR^d$} &
                                                     points in the
  Lie quadric \\
\hline
  \hline
 \multicolumn{3}{|c|}{{\bf Incidence}} \\
  \hline
1 &  Point $\bfp$ is incident to the oriented hypersphere $S$ &
                                                         \multirow{2}{*}{$B(\xi', \sigma')
                                                         = \la \xi, \sigma \ra = 0$} \\
                                \hline
2&   Point $\bfp$ is incident to the oriented hyperplane $H$ & \multirow{2}{*}{$B(\xi', \pi')
                                                         = \la \xi, \pi \ra = 0$} \\
                                \hline
  \multicolumn{3}{|c|}{{\bf Inclusion} } \\
 \hline
3 &  Point $\bfp$ is  inside the unoriented
  hypersphere $S$ with $r > 0$  &
                                                     \multirow{2}{*}{$B(\xi', \sigma')
                                                         = \la \xi, \sigma \ra > 0$} \\
                                \hline
4 &   Point $\bfp$ is  outside the unoriented
  hypersphere $S$ with $r > 0$  &
                                                     \multirow{2}{*}{$B(\xi', \sigma')
                                                         = \la \xi, \sigma \ra < 0$} \\
                                \hline
 5 &  Point $\bfp$ is the positive
  half-space $H^+$ defined by oriented hyperplane  $H$ &  \multirow{2}{*}{$B(\xi', \pi')
                                                         = \la \xi, \pi \ra > 0$} \\
  \hline
 6& Sphere $S$ is a subset of the  positive half-space $H^+$ defined by oriented hyperplane  $H$ & \multirow{3}{*}{$\la \sigma, \pi \ra > 0$} \\
  \hline
 7 & Sphere $S(\bfq_1,r_1)$ is a subset of  the exterior of sphere $S(\bfq_2,r_2)$ & \multirow{2}{*}{$\la \sigma_{\bfq_1,r_1}, \sigma_{\bfq_2,-r_2} \ra < 0$} \\
  \hline
  \multicolumn{3}{|c|}{{\bf M\"obius scalar product} }\\
  \hline
8 &   M\"obius scalar product $\rho(S_1,S_2) $ & \multirow{2}{*}{$\la \sigma_1, \pi(\sigma_2)\ra \le 0$} \\
& is nonpositive & \\
  \hline
  \multicolumn{3}{|c|}{{\bf Tangency and oriented contact} }\\
  \hline
9 &   Oriented hyperspheres $S_1,S_2$ are in oriented contact & \multirow{2}{*}{$\la \sigma_1, \sigma_2 \ra = 0$} \\ 
  \hline
  10 & Oriented hypersphere $S$  and oriented hyperplane $H$ are in oriented contact  
  &  \multirow{2}{*}{$\la \sigma, \pi
    \ra = 0$} \\
  \hline
11 &   Unoriented hyperspheres $S(\bfq_1,r_1)$  and   $S(\bfq_2,r_2)$ are
   externally tangent 
  &  \multirow{2}{*}{$\la \sigma_{\bfq_1,r_1}, \sigma_{\bfq_2,-r_2}
  \ra = 0$} \\
  \hline
\end{tabular}
\smallskip
    \caption{Translation of properties of subsets in 
      $\boldR^d$ to properties of corresponding points in the Lie
quadric $\calQ^{d+1}$ in standard form, in terms of the Lie inner
product and Minkowski inner product.  The point $\bfp \in \boldR^d$
corresponds to $[\xi] \in \calQ^{d+1},$ the oriented hyperspheres $S,
S_1, S_2 \subseteq \boldR^d$ correspond to $[\sigma], [\sigma_1]$ and
$[\sigma_2]$ in $\calQ^{d+1}$, and oriented hyperplanes $H, H_1, H_2
\subseteq \boldR^d$ correspond to $[\pi], [\pi_1]$ and $[\pi_2]$ in
$\calQ^{d+1}$.  The projections of $\xi, \sigma$ and $\pi$ to the first
$d+2$ coordinates in $\boldR^{d+2}$ are given by $\xi', \sigma'$ and
$\pi'$ respectively.}\label{table: Lie properties}
    \end{table}

 \begin{table}   \setlength{\extrarowheight}{2pt}
\centering
\begin{tabular}{| l | c |}
  \hline
  \hfill Distance in $\boldR^d$ \hfill & Lie inner product in $\boldR^{d+3}$ \\
  \hline
  \hline
  $\dist(\bfp,\bfq)$ & $-2 \la \xi_p, \xi_{\bfq}\ra$ \\
 \hline
  $\dist(\bfp,H_{\bfn,h})$ for $\bfp \in H^+_{\bfn,h}$ & $\la \xi_p,
                                                         \pi_{\bfn,h}\ra$
  \\
    \hline
  $\dist(S(\bfq,r),H_{\bfn,h})$ for $S(\bfq,r) \subseteq H^+_{\bfn,h}$
                     & $\la \sigma_{\bfq,r}, \pi_{\bfn,h}\ra$ \\
    \hline
  $\Pow(\bfp, S(\bfq,r))$ & $-2\la \sigma_{\bfq,r}, \xi_{\bfp}\ra$ \\
  \hline 
\end{tabular}
\smallskip
    \caption{Distances between objects in $\boldR^d$ and 
      the power of a point with respect to a sphere expressed in terms
      of the Lie inner product and Lie
      spheres in standard  coordinates.
    }
      \label{table: distances}
    \end{table}

    \subsection{Connections between M\"obius geometry and Lie sphere
      geometry}
    The M\"obius sphere $\Sigma^d$ is embedded in the Lie quadric
$\calQ^{d+1}$ as the set of points with signed radius zero:
 \[ \Sigma^d \cong [\bfe_{d+3}]^\perp \cap \calQ^{d+1} = \{ [\bfx] \in
\calQ^{d+1} \, : \, \la \bfx, \bfe_{d+3}\ra = 0 \} .  \] 

 A condition on a point $\bfx$ in $\boldR^d$ expressed in terms of
$\xi_{\bfx}'$ in $\boldR^{d+2}$ and the Lorentz inner product
translates to an analogous condition on $\xi_{\bfx}$ in $\boldR^{d+3}$
in terms of the Lie inner product because the last coordinate of
$\xi_\bfx$ is zero.  For any oriented sphere $\sigma_{\bfq,r}$ with
corresponding unoriented sphere $\sigma_{\bfq,r}'$, $\la \xi_\bfx,
\sigma_{\bfq,r}\ra = B(\xi_\bfx', \sigma_{\bfq,r}') $ and $ \la
\xi_\bfx, \pi_{\bfn,h}\ra = B(\xi_\bfx', \pi_{\bfn,h}'),$ where
$\pi_{\bfn,h}$ is an oriented hyperplane and $\pi_{\bfn,h}'$ is same
hyperplane without orientation.

For any $\bfx = (x_i) \in \boldR^{d+3}$ let
\[ \pi(\bfx) = (x_1, x_2, \ldots, x_{d+2}, 0) \] denote its orthogonal
projection to the coordinate plane $x_{d+3}=0$.  For the standard Lie
coordinates $\sigma_{\bfq,r}$ of a Lie sphere $S(\bfq,r)$ in
$\boldR^d$,
\begin{equation}\label{eqn: minkowski to lie}
  \pi(\sigma_{\bfq,r}) = \left(\frac{1 +\bfq \cdot
        \bfq -  r^2}{2},\frac{1 - \bfq  \cdot \bfq + r^2}{2}, \bfq,
      0 \right).  \end{equation}
The Lorentz inner product of unoriented spheres $\sigma_{\bfq_1,r_1}'$
and $\sigma_{\bfq,r}'$ in standard coordinates in $\boldR^{d+2}$ is
equal to the Lie inner product of  projections of
$\sigma_{\bfq_1,r_1}$ and $\sigma_{\bfq,r}$ in $\boldR^{d+3}$:
\begin{equation}\label{eqn: minkowski to lie inner product}
 B(\sigma_{\bfq_1,r_1}', \sigma_{\bfq_2, r_2}') = \la \pi(\sigma_{\bfq_1,r_1}),
 \pi(\sigma_{\bfq,r})\ra.\end{equation}

The M\"obius scalar product of spheres $S(\bfp_1,r_1)$ and
$S(\bfp_2,r_2)$ as in Equation \eqref{eqn: Mobius} is simply the
Lorentzian inner product
\[ \rho(S(\bfp_1,r_1), S(\bfp_2,r_2)) = B(\sigma_{\bfp_1,r_1}',
  \sigma_{\bfp_2,r_2}')\]
of the corresponding points $\sigma_{\bfp_1,r_1}'$ and
$\sigma_{\bfp_2,r_2}'$ in $\boldR^{d+2}$. But $B(\sigma_{\bfp_1,r_1}',
\sigma_{\bfp_2,r_2}') = \la \pi(\sigma_{\bfp_1,r_1}),
\sigma_{\bfp_1,r_1}\ra$.  Thus the M\"obius scalar product may be
expressed in terms of the Lie inner product as
\[  \rho(S(\bfp_1,r_1), S(\bfp_2,r_2)) =  \la \pi(\sigma_{\bfp_1,r_1}),
\sigma_{\bfp_1,r_1}\ra. \]
\subsection{The space of spheres, in terms of standard Lie
  coordinates}\label{space of spheres}
We would like to describe the space of spheres as a subset of
Euclidean space rather than projective space.  To do that, we use
standard coordinates.  From Table \ref{table: Lie coordinates
properties}, we see that the set $\Sph(\boldR^{d})$ of all
representatives for oriented spheres and point spheres in standard
form is
\[ \Sph(\boldR^{d}) = \{ \bfx \in \boldR^{d+3} \, : \, x_1 + x_2 =
  1, \la \bfx, \bfx \ra = 0 \}.\]
We have already seen that the paraboloid
$\calR$ in Equation \eqref{eqn: paraboloid}
parametrizes the set of point spheres in $\Sph(\boldR^d)$. 
On the other hand, we can identify the set of oriented spheres with
$\boldR^d \times \boldR$ in the obvious way, assigning $S(\bfq,r)$ to
$(\bfq,r)$.  Equation \eqref{eqn: lie sphere} defines how to map
$(\bfq,r)$ to $\bfsigma_{\bfq,r} = (\sigma_i)$ in $\Sph(\boldR^{d}) $:
\begin{align*}
  \sigma_1
  &=    \smallfrac{1}{2}(1 +\| \bfq \|^2 -  r^2),  \\
 \sigma_2&=
           \smallfrac{1}{2}(1 - \| \bfq \|^2 + r^2) \\
  (\sigma_3, \ldots,\sigma_{d+2}, \sigma_{d+3}) &= (\bfq, r).
\end{align*}
Conversely, to map $\bfsigma_{\bfq,r} = (\sigma_i)$ in to $(\bfq,r)$,
one simply projects to the coordinates $ (\sigma_3,
\ldots,\sigma_{d+2}, \sigma_{d+3}) .$ Define the {\em sphere center
map} $\prc: \boldR^{d+3} \to \boldR^d$ and {\em radius map} $r:
\boldR^{d+3} \to \boldR^d$, by
\begin{equation}\label{eqn: prc}
  \prc:   (x_1, x_2, \ldots, x_{d+3})      \mapsto  (x_3, \ldots, x_{d+2}).
      \end{equation}
and
\begin{equation}\label{eqn: radius}
  r:  (x_1,x_2,x_3, \ldots, x_{d+2}, x_{d+3})  \mapsto
      x_{d+3}. \end{equation} 
Note that for $\bfsigma_{\bfq,r} = (\sigma_i)$ in $\Sph(\boldR^{d}) ,$
the radius is given by
    \[ \sigma_{d+3} = \pm \sqrt{-\sigma_1^2 + \sigma_2^2 + \cdots +
\sigma_{d+2}^2},\] where the sign choice depends on the sign of the
radius in $\bfsigma_{\bfq,r}.$

\section{Generalized Voronoi diagrams defined by extremal 
  spheres}\label{section: extremal spheres}
In this section, we prove Theorem \ref{thm: empty sphere} and give
applications of the theorem.

\subsection{Proof  of Theorem \ref{thm: empty sphere} }

Having set out the basics of Lie sphere geometry, the proof of Theorem
\ref{thm: empty sphere} will follow relatively easily.  Let $\calP$ be
a finite set of points in $\boldR^d,$ let $\calH$ be a set of
half-spaces in $\boldR^d$, and let $\calS$ and $\calT$ be finite sets
of spheres in $\boldR^d,$ where not all of these sets are empty.  Assume that
$\Sph(\calD)$ is nonempty. Let
$Z(\calD)$ be the set of sphere centers in $\boldR^d$ determined by
the data set $\calD=(\calP^-,\calP^+, \calH, \calS, \calT)$ as
described in Definition \ref{defn: data set}. 

We present the set of inequalities in $\boldR^{d+3}$ determined by the
data set.  First we convert all the sites into Lie sphere coordinates:
for a point $\bfx$ in $\calP^- \cup \calP^+$, let $\xi_{\bfx}$ be its
standard coordinates; for a half-space $H_{\bfn, h}^+$, associate the
oriented hyperplane $\pi_{\bfn,h}$; to an unoriented sphere
$S(\bfq,t)$ in $\calS$ associate $\sigma_{\bfq,t}$ (which has positive
orientation); and to an unoriented sphere $S(\bfq,t)$ in $\calT$,
associate the sphere $\sigma_{\bfq,-t}$ (which has negative
orientation).  For an arbitrary unoriented sphere $S(\bfx,r)$ with
center $\bfx$ and radius $r$, endow it with positive orientation and
use $\sigma_{\bfx,r}$ to denote its standard coordinates in
$\Sph(\boldR^d)$.  Next, for each point representing a site in Lie
sphere coordinates, define a non-strict linear inequality in $(d+3)$
variables $\sigma = (\sigma_i)$ expressed in terms of the Lie inner
product as laid out in Table \ref{table: inequalities}.

\begin{example} Let $\calD=(\calP^-,\calP^+, \calH, \calS, \calT)$ be
a data set.  Suppose that the sphere $S( (3,1), 5) \subseteq \boldR^2$
is in in the set $\calS.$ The Lie sphere coordinates associated to $S(
(3,1), 5) $ are $\sigma_{(3,1),5} = (-7, 8, 3,1,5)$, and the
inequality associated to $S( (3,1), 5) $ from Table \ref{table:
inequalities} is
 \[ \la  (-7, 8,
   3,1,5), (\sigma_1,\sigma_2,\sigma_3,\sigma_4,\sigma_5) \ra
   = 7 \sigma_1 +  8 \sigma_2 + 
  3\sigma_3 + \sigma_4 - 5\sigma_5 \le 0.\]
\end{example}

\begin{table}
    \setlength{\extrarowheight}{2pt}
\centering
\begin{tabular}{| l | c |}
  \hline
  Site defining a geometric& Associated inequality   \\
  condition on $S(\bfx,r)$  & in variable $\sigma= \sigma_{\bfx,r}$\\
  \hline
  \hline
  Point $\bfp$ in $\calP^-$  \qquad & \quad $\la \sigma, \xi_{\bfp} \ra
                                        \le  0$ \quad \\
  \hline
  Point $\bfp$ in $\calP^+$ & $\la \sigma, \xi_{\bfp} \ra \ge  0$ \\
  \hline
   Half-space $H^+_{\bfn, h}$ in $\calH$& $\la \sigma, \pi_{\bfn, h} \ra \le  0 $\\
  \hline
  Sphere $S_{\bfq,t}$ in $\calS$ & $\la \sigma, \pi(\sigma_{\bfq,t}) \ra 
     \le  0 $\\
  \hline
  Sphere $S_{\bfq,t}$ in $\calT$ &$\la \sigma, \sigma_{\bfq,-t} \ra \le  0 $ \\
  \hline
\end{tabular}
\smallskip
    \caption{The system of inequalities encoding the geometric  conditions in
      Definition \ref{defn: data set}.}
    \label{table: inequalities}
    \end{table}

Now we are ready to prove the theorem. 
\begin{proof}  [Proof of Theorem \ref{thm: empty sphere}]
 Let $\Sph(\calD)$ be the set of points in $\boldR^d$ determined by
the data set $\calD=(\calP^-,\calP^+, \calH, \calS, \calT)$ as
described in Definition \ref{defn: data set}.  Let $S(\bfx,r)$ be an
arbitary unoriented sphere in $\boldR^d$ in the set $\Sph(\calD)$,  and let $\sigma=\sigma_{\bfx,r}.$ 
For each site in the data set, associate the linear inequality as in Table \ref{table: inequalities}.  Let $P$ denote the polyhedron
that is the intersection of the half-spaces defined by the
inequalities indexed by the data set.

Now we apply facts from Table \ref{table: Lie properties} to the
sphere $S(\bfx,r)$:
\begin{enumerate}
\item For all points $\bfp$ in $\calP^-$, $\bfp$ is not in the
interior of $S(\bfx,r)$ if and only if $\la \sigma, \xi_{\bfp} \ra \le
0$.  (By Row 3 of Table \ref{table: Lie properties}.)
\item  For all points $\bfp$ in $\calP^+,$ $\bfp$ is not in the
exterior of $S(\bfx,r)$ if and only if $\la \sigma, \xi_{\bfp} \ra \ge
0$. (By Row 4 of Table \ref{table: Lie properties}.)
\item  For all half-spaces $H^+_{\bfn,h}$ in $\calH$, the sphere
$S(\bfx,r)$ is a subset of the closed half-space
$\overline{H^+_{\bfn,h}}$ if and only if $\la \sigma, \pi_{\bfn,h} \ra
\le 0$. (By Rows 6 and 10 of Table \ref{table: Lie properties}.)
\item  For all spheres $S(\bfq,t)$ in $\calS$, the M\"obius scalar
product of $S(\bfx,r)$ and $S(\bfq,t)$ satisfies
     \[ \rho(S(\bfx,r),S(\bfq,t)) \le 0\] if and only if $\la \sigma,
\pi( \sigma_{\bfq,t}) \ra \le 0$.  (By Row 8 of Table \ref{table: Lie
  properties}.)
\item  For all $S(\bfq,t)$ in $\calT $, the sphere $S(\bfx,r)$ is a
subset of the closure of the exterior of $S(\bfq,t)$ if and only if
$\la \sigma, \pi(\sigma_{\bfq,-t}) \ra \le 0$. (By Rows 7 and 11 of
Table \ref{table: Lie properties}.)
\end{enumerate}
It follows that that the Lie sphere $\sigma$ is in the polyhedron $P$
if and only if the sphere $S(\bfx,r)$ satisfies
the associated geometric conditions.  Because $\Sph (\calD)$ is
nonempty, $P$ is nonempty.  This proves the first two parts of the
theorem.

The third part of the theorem statement describes extremal spheres.
To see that a sphere $S(\bfx,r)$ is extremal occurs for a condition if
and only if equality holds in the corresponding inequality, we again
refer to Table \ref{table: Lie properties}.
\begin{enumerate}
\item The condition for a point $\bfp$ in $\calP^+$ is extremal if
$\bfp$ is incident to $S(\bfx,r)$.  This holds if and only if $\la
\sigma, \xi_{\bfp} \ra = 0$. (By Row 1 of Table \ref{table: Lie
  properties}.)
\item The condition for a point $\bfp$ in $\calP^-$ is extremal if
$\bfp$ is incident to $S(\bfx,r)$.  This holds if and only if $\la
\sigma, \xi_{\bfp} \ra = 0$. (By Row 1 of Table \ref{table: Lie
properties}.)
\item The condition for a half-space $H^+_{\bfn,h}$ in $\calH$ is
extremal if and only if the sphere $S(\bfx,r)$ is tangent to the
boundary of the closed half-space $\overline{H^+_{\bfn,h}}$.  This
occurs if and only if $S(\bfx,r)$ is in oriented contact with the
oriented hyperplane.  This occurs if and only if $\la \sigma,
\pi_{\bfn,h} \ra = 0$. (By Row 10 of Table \ref{table: Lie
  properties}.)
\item The condition for a sphere $S(\bfq,t)$ in $\calS$ holds
extremally if and only if
\[ 0 = \rho(S(\bfx,r),S(\bfq,t)) = \la \sigma_{\bfx,r},
     \pi(\sigma_{\bfq,t} \ra.\]
\item  The condition for a sphere $S(\bfq,r)$ in $\calT$ holds
extremally if and only if $S(\bfq,r)$ and $S(\bfx,r)$ are in
externally tangent, in which case $ \la \sigma_{\bfx,r},
\pi(\sigma_{\bfq,-t} \ra = 0$. (By Row 11 of Table \ref{table: Lie
properties}.)
\end{enumerate}

The fourth part of the theorem statement follows from the fact that the
map $\prc$ projects $\sigma_{\bfx,r} \in \Sph(\boldR^d)$ to the
center  $\bfx$ in $\boldR^d.$ 
\end{proof}

\subsection{Applications of Theorem \ref{thm: empty sphere}}
In this section, we present direct applications of Theorem
\ref{thm: empty sphere}, starting with the classical Voronoi diagram
for a set of points in $\boldR^d$.

\begin{cor}\label{cor: classical Voronoi}
Let $\calP = \{\bfp_1, \bfp_2, \ldots, \bfp_n\}$ be a set of 
points in $\boldR^d.$ The classical Voronoi diagram for $\calP$ is
obtained from Theorem \ref{thm: empty sphere} using the data set
$(\calP^-,\calP^+, \calH, \calS, \calT)$, where $\calP^- =
\calP$ and the remaining sets $\calP^+, \calH, \calS$ and $\calT$
are empty.  The corresponding system of inequalities in $\boldR^{d+3}$
is
  \[  \la \xi_{\bfp_i} , \sigma \ra \ge 0 \qquad \text{for $i=1, \ldots, n$}.\]
\end{cor}

\begin{proof} The classical Voronoi diagram is the minimization
diagram for the distance functions $\dist(\bfp_i, \cdot)$.  For a
subset $I$ of $[n]$ a point $\bfx$ is in $\Cell(I)$ if and only if
$\bfx$ is the center of a sphere $S(\bfx,r)$ incident to $\bfp_i$ for
$i \in I$, with all other points outside the sphere $S(\bfx,r).$
\end{proof}

In the next corollary we describe the farthest point generalized
Voronoi diagram.

\begin{figure}
  \begin{center}
  \includegraphics[scale = .6]{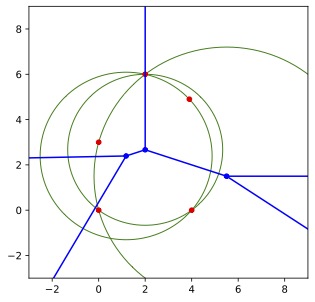} 
\end{center}
\caption{The farthest point Voronoi diagram}\label{fig: farthest point}
\end{figure}

\begin{cor}\label{cor: farthest point}
Let $\calP = \{\bfp_1, \bfp_2, \ldots, \bfp_n\}$ be a set of points in
$\boldR^d.$ The farthest point Voronoi diagram for $\calP$ is
described by Theorem \ref{thm: empty sphere} using the data set
$(\calP^-,\calP^+, \calH, \calS, \calT)$, where $ \calP^+ = \calP$ and
the remaining data sets $\calP^-, \calH, \calS$ and $\calT$ are empty.
The corresponding system of inequalities is inequalities
\[  \la \xi_i , \sigma \ra \le 0 \qquad \text{ for $i=1, \ldots, n$}.\]
\end{cor}
See Figure \ref{fig: farthest point} for an example.

\begin{proof}
Let $\bfx$ be a point in $\boldR^d$.  The point $\bfp_k$ is the farthest
point in $\calP$ from $\bfx$ if and only if there exists a sphere
$S(x,r)$ of radius $r$ centered at $x$ such that $\bfp_k$ is incident to
the sphere and for all $i$, $p_i$ is in $S(x,r)$ or incident to it.
\end{proof}
Now we consider the case of spherical sites.  Let $\calT= \{ S(\bfq_i, r_i)\}$ be
a finite set of spheres in $\boldR^d$ with positive radii.  The {\em
  Apollonius diagram}  is the minimization diagram $\Vor(\calT)$  for the distance
functions $f_i(x) = \dist(\bfx, S_i(\bfq_i,r_i))$.  We take $\Dom(\calT)=\boldR^d
\setminus \cup_i  \setint(S(\bfq_i, r_i)) $ as the common domain of the
functions $f_i$.  (We measure the distance from a point $\bfp$in $\boldR^d$ to a
compact subset $K$ of $\boldR^d$ with the function $\dist(x,K) = \min \{ \dist(\bfx,\bfy) \,
: \, \bfy \in K\}.$) 
  
\begin{figure}
  \begin{center}
  \includegraphics[scale = .6]{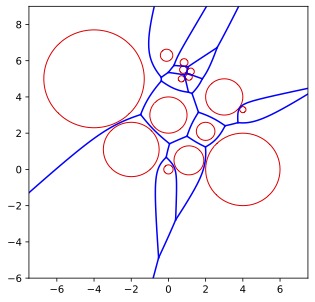} 
\end{center}
\caption{The Apollonius diagram.}\label{fig: apollonius}
\end{figure}

\begin{cor}\label{cor: spherical sites}
  Let $\calT_0= \{S_1, S_2, \ldots, S_n\}$ be a set of spheres in $\boldR^d.$
 The
Voronoi diagram $\Vor(\calT_0)$  is described by Theorem \ref{thm:
empty sphere} using the data sets $(\calP^-,\calP^+, \calH, \calS,
\calT)$, with $ \calT = \calT_0$, and with the remaining data sets  empty.
The corresponding system of inequalities is
\[ \la \sigma_i , \sigma \ra \le 0 \text{ for $i=1, \ldots, n$}.\]
\end{cor}

\begin{proof}
  Let $\bfx$ be in the domain $\Dom(\calT)$.
  For any $i,$ let $f_i (\bfx)= 
= \dist( \bfx, S_i(\bfq_i,r_i)) $.  Fix $j \in [n].$  The function
$f_j$ is minimized by $\bfx$ if and only if among all the spheres
$S(\bfx,r) \subseteq \Dom(\calT)$ , there is a
sphere centered at $x$ externally tangent to $S(\bfq_j, r_j)$.  This is
equivalent to Condition (5) holding for all $i$ and holding extremally
for the index $j$; or, in Lie sphere coordinates,
\[ \la \sigma_i , \sigma \ra \le 0 \text{ for $i=1, \ldots, n$},\]
with equality for the equation indexed by $j.$
\end{proof}

\begin{figure}
  \begin{center}
  \includegraphics[scale = .6]{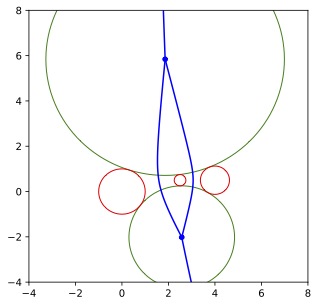} 
\end{center}
\caption{The Voronoi diagram for three small circles (in red).  Here the
  minimization region for the
  index set   $\{1,2,3\}$  is the disjoint union of two points,
  reflecting the fact that the
  intersection of the corresponding subspace in $\boldR^5$ intersects the Lie quadric
  in two distinct points.}\label{fig: 2 point intersection}
\end{figure}

When sites are spheres,  there may be
multiple edges between two vertices of the generalized Voronoi diagram.  Note $\Cell({1,2,3})$
has  is a set of cardinality two in Figure \ref{fig: 2 point intersection}.
  
\begin{cor}[\cite{aurenhammer1987power}]
  \label{cor: power diagram 1}
Let $\calS_0 = \{ S(\bfq_i, r_i) \}$ be a set of $n$ spheres in
$\boldR^d$ whose interiors are mutually disjoint.  The power diagram for $\calS_0$ is described by
Theorem \ref{thm: empty sphere} using the data set $(\calP^-,\calP^+,
\calH, \calS, \calT)$, where $\calS = \calS_0$ and the remaining sets $\calP^-, \calP^+,
\calH$ and $\calT$ are empty.
The corresponding system of inequalities is
  \[ \la \sigma_{\bfq_i,r_i}' , \sigma \ra \le 0 \text{ for $i=1,
      \ldots, n$},\]
where $\sigma_{\bfq_i,r_i}$ is the Lie sphere coordinates for the sphere $S_i$,
in standard form.
\end{cor}

\begin{figure}
  \begin{center}
  \includegraphics[scale = .6]{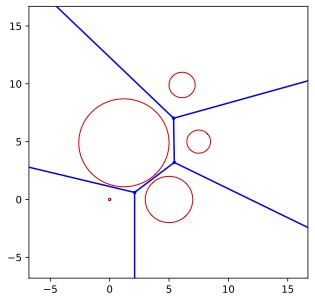} 
\end{center}
\caption{Power diagram}\label{fig: power diagram}
\end{figure}

\begin{figure}
  \begin{center}
  \includegraphics[scale = .6]{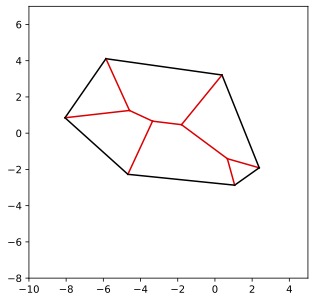} 
\end{center}
\caption{The medial axis of a convex polygon}\label{fig: medial axis}
\end{figure}

We will view the power diagram in terms
 of centers of maximally mostly hollow spheres, where the sphere
 $S(\bfp,r)$ is considered {\em mostly hollow} relative to
 $S(\bfq,t)$ if $\| \bfp - \bfq \|^2 -r^2 -t^2 \ge 0,$ or
 equivalently, $(\sigma_{\bfp,r},\sigma_{\bfq,t}) \le 0.$ This
 occurs if and only if 
\[  \la \sigma_{\bfp,r}',\sigma_{\bfq,t}\ra = (\sigma_{\bfp,r},\sigma_{\bfq,t}) =
  -\smallfrac{1}{2} (\Pow(\bfp, S(\bfq,t)) - r^2) \le 0\]
 Thus, the
 sphere  $S(p,r)$ is mostly hollow relative to $S(\bfp_1, r_1) $ if
 and only if $\Pow(\bfp, S(\bfp_1,r_1)) \ge r^2$.
\begin{proof}
Let $\calS_0 = \{ S(\bfq_i, r_i) \}$ be a set of spheres as in the
statement of the theorem.  Let $\bfp$ be a point so that
\[ R_\bfp := \Pow(\bfp, S(\bfp_j,r_j)) = \min_i \Pow(\bfp,
  S(\bfp_i,r_i)). \]
Because $\bfp \not \in \cup_i \setint S(\bfp_i,r_i)$, we know that
$R_\bfp \ge 0.$
Then $S(\bfp,R_\bfp)$ is a maximally hollow sphere about $\bfp$ and$\la
\sigma_{\bfp_j,r_j}',\sigma_{\bfp,r_\bfp}\ra  \le 0.$
\end{proof}
Next we apply Theorem \ref{thm: empty sphere}  to the medial  axis of
a  convex polygon.   
\begin{cor}\label{cor: medial axis}
Let $T = \cup H_{\bfn_i,h_i}^+$ be a polyhedron defined as the
intersection of half-spaces $H_{\bfn_i,h_i}^+,$ for $i \in [n]$.  Let
$\pi_{\bfn_i,h_i} $ be the standard Lie coordinates for the half-space
$H_{\bfn_i,h_i}^+.$ The medial axis for $T$ is defined by Theorem
\ref{thm: empty sphere} using the data set $(\calP^-,\calP^+, \calH,
\calS, \calT)$, where $\calP^-, \calP^+ , \calS $ and $\calT$ are
empty and $\calH = \{ \pi_{\bfn,h}\}$ The corresponding system of
inequalities is inequalities is
\[  \la \xi_i , \pi_{\bfn, h} \ra \le 0 \text{ for $i=1, \ldots, n$}.\]
\end{cor}
See Figure \ref{fig: medial axis}.

\begin{proof}  A point $\bfx$ in $\boldR^d$ is in the medial axis for
$T$ if and only if there is a hypersphere (outwardly oriented)
centered at $\bfx$ which is tangent to two or more of the oriented
hyperplanes $H_{\bfn_i,h_i}$.  
\end{proof}

\begin{figure} 
  \begin{center}
      \includegraphics[scale = .6]{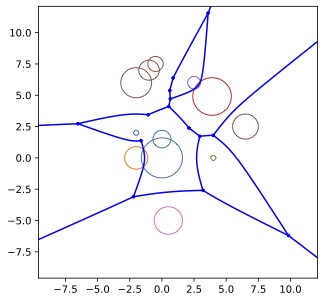} 
\end{center}
\caption{Mixed sites, as described in Example \ref{example: mixed sites}.}\label{fig: mixed}
\end{figure}

\begin{example}\label{example:  mixed sites}
  In Figure \ref{fig: mixed} we show an example where the sites in
  $\calD$ are of mixed type.  In fact, we have allowed some sites to
  be unions of objects,  and removed edges accordingly to obtain the
  minimization diagram. Objects of the same color are considered to be
  part of the same generalized site. 
\end{example}

\section{Minimization diagrams from restrictions of affine
  functions}\label{section: affine functions}

\subsection{Observations about minimization diagrams for compositions
  of affine functions and embeddings.}

Before we state the next theorem we make some simple observations
about minimization diagrams involving embeddings and intersections.

Suppose that $\{L_i\}$ is a family of affine functions from $\boldR^n$
to $\boldR.$ The minimization diagram for the functions in $\calL =
\{L_i\}$ is straightforward to compute as the functions are affine;
each cell is the intersection of half-spaces.  Suppose that $N$ is a
closed subset of $\boldR^n.$ For each affine function $L_i,$ let
$L_i|_N: N \to \boldR$ denote its restriction to $N$, and let
$\calL^N=\{L_i|_N\}$. Each minimization region for the minimization
diagram $\MD(\calL^N)$ is the intersection of $N$ with the
corresponding minimization region for $\MD(\calL)$: The minimization
region $\MRN(I)$ in $\MD(\calL^N)$ for an index set $I$ as defined in
Equation \eqref{defn MR} is equal to the intersection of $N$ with the
corresponding polyhedral cell in $\MD(\calL)$.  To be precise,
\begin{equation}\label{eqn: MR int} \MD(\calL^N) = \{ \overline{\MRN(I)}  \}  =  \{
  N \cap \overline{\MRI } \, :
 \overline{\MRI} \in  \MD(\calL) \},
\end{equation}
where $I$ varies over all possible sets of indices.  Thus, the
minimization diagram $\MD(\calL^N)$ for the family of functions
$\calL^{N} = \{L_i|_N\}$ may be obtained intersecting the polyhedra in
$\MD(\calL) $ with $N.$

\begin{figure}\label{fig: cd}
  \begin{center}
    \hfill
\begin{tikzcd}
 \boldR^d \arrow[r, "\phi"] \arrow[dr, "f_i"']
& \boldR^{n} \arrow[d, "L_i"]\\
& \boldR
\end{tikzcd} \hfill \begin{tikzcd}
\boldR^d \arrow[r, "\phi"] \arrow[dr, "f^-=\min_i \{f_i\}"']
& \boldR^{d+3} \arrow[d, "L^-=\min_i \{L_i\}"]\\
& \boldR
\end{tikzcd}
\hfill
\strut
\end{center}
\caption{Functions $f_i$ factor through affine maps $L_i$}
\end{figure}
\hfill
\strut

Next, suppose that there is an embedding $\phi: \boldR^d \to
\boldR^n$, where $n > d$, and suppose that the image $N =
\phi(\boldR^d)$ is closed.  Let $\calF = \{ f_i \}$ be a set of
real-valued functions with common domain $\boldR^d$ such that each one
is the composition $f_i = L_i \circ \phi$, where $L_i: \boldR^n \to
\boldR$ is an affine function.  See Figure \ref{fig: cd}. In this
situation, for all $\bfx$ in $\boldR^d,$
\[ \min_i \{ f_i (\bfx) \} = \min_i \{ (L_i \circ \phi ) (\bfx )\} =
\min_i \{ L_i (\bfy)\},\]
where $\bfy = \phi(\bfx)$. As a consequence, the sets of minimizing
indices coincide; in the notation as in Equation \eqref{defn I}, the
set $I(\bfx)$ for $\calF$ equals the set $I(\phi(\bfx))$ for $\calL =
\{ L_i \}$. Hence, $\bfx$ is in the minimization region $\MRF(I)$ for
$\{ f_i\}$ if and only if $\phi(\bfx)$ is in the minimization region
$\MRN (I)$ for $\{ L_i|_N\}.$ It follows that for any index set $I$,
the minimization region $\MRF (I)$ is equal to the preimage
$\phi^{-1}({\MRN}(I))$ of the corresponding minimization region in
$N$.  We rewrite using Equation \eqref{eqn: MR int} to obtain
\[ \MRF(I) = \phi^{-1}(\MRN(I)) = \phi^{-1} (N \cap \overline{\MR(I) } ). \]
Note that $\overline{\MRF(I)} = \overline{\phi^{-1}(\MRN(I))} =
\phi^{-1}(\overline{\MRN(I)} ) $ as $\phi$ was assumed to be an
embedding.  We conclude that the miniization diagram for $\calF$ is
\[  \MD(\calF) = \{ \overline{\phi^{-1}(\MRN(I) )} \}  =  \{
  \phi^{-1}(N \cap \overline{\MRI } )\, :
 \overline{\MRI} \in  \MD(\calL) \},
\]
where $I$ varies over all possible sets of indices.  The upshot of
this discussion is that the minimization diagram for $\calF$ may be
computed by first finding the polyhedral minimization regions for
$\MD(\calL)$, and then for each of these regions, intersecting with
$N$ and taking the pre-image under $\phi.$ The following proposition
states this precisely.

\begin{prop}\label{prop: MD restriction}
Let $\phi: \boldR^d \to \boldR^n$ be an embedding of $\boldR^d$ into
$\boldR^n$ with closed image $N$.  Let $\calF = \{ f_i\}_{i=1}^n$ be a
family of real-valued functions  with domain $\boldR^N$. Suppose
that for each $i$ there is an affine function $L_i: \boldR^n \to \boldR$
so that $f_i = L_i \circ \phi.$ The minimimization diagram for $\calF$
is the set of minimization regions
\[  \MD(\calF) = \{ \overline{\phi^{-1}(\MRN(I) )} \}  =  \{
  \phi^{-1}(N \cap \overline{\MRI } )\, :
 \overline{\MRI} \in  \MD(\calL) \}.
\] 
\end{prop}

Now we are set to prove Theorem \ref{thm: affine on quadric}. 
\begin{proof}[Proof of Theorem \ref{thm: affine on quadric}]
  Let $\calF = \{f_i\}_{i=1}^n$ and $\calL = \{ L_i \}_{i=1}^n$ be as
in the statement of the theorem. Let $N=Q^d.$
   \end{proof}

Next we prove Proposition \ref{prop: of form}. 
\begin{proof}[Proof of Proposition \ref{prop: of form}]
Let $f: \boldR^d \to \boldR$ be a function of form 
     \[  f(\bfx) =  a \| \bfx - \bfq \|^2 + \bfb     \cdot \bfx + c,
       \quad \text{for $\bfx \in \boldR^d$}, \]
where $a, c \in \boldR$ and $\bfb \in \boldR^d$.  We rewrite $f(\bfx)$
as
     \begin{align}
       f(\bfx) &= a \| \bfx \|^2 - 2a \, \bfq \cdot \bfx + a \| \bfq\|^2 +
                 \bfb     \cdot \bfx + c\\ \notag
       &= a \| \bfx \|^2 + ( \bfb - 2a \bfq ) \cdot \bfx + a \|
         \bfq\|^2 + c . \label{eqn: coeffs}
     \end{align}
We  make a linear change of variables from parameters $(a, \bfb, c)$
to $\bfa = (a_i)$ in $\boldR^{d+3}$ by letting

\begin{equation}\label{eqn: linear cov}
\begin{split}
  \begin{bmatrix}
    a_1 \\ a_2 
  \end{bmatrix}  &=
  \begin{bmatrix}
    -1 & -1  \\
    -1 & 1
  \end{bmatrix} 
  \begin{bmatrix}
    a \\ a \| \bfq\|^2 + c 
  \end{bmatrix},  \\
  (a_3, \ldots, a_{d+2} ) &= \bfb - 2a  \bfq, \enskip \text{and}  \\
  a_{d+3}&=0.
           \end{split}
  \end{equation}
  
This allows us to rewrite $f$ as 
\begin{equation}\label{eqn: f}
           f(\bfx) =
                     \smallfrac{1}{2}(a_2-a_1)-\smallfrac{1}{2}(a_2+a_1)
                     \| \bfx \|^2 +  (a_3, \ldots, a_{d+2} )  \cdot
                     \bfx.
\end{equation}
Recall that the image of a point $\bfx$ in $\calR$
under $\phi$ is
  \[ \phi(\bfx) = \left( \frac{1+\bfx \cdot \bfx}{2}, \frac{1 - \bfx
        \cdot \bfx}{2}, \bfx, 0\right). \]
 A short computation shows that $ f(\bfx) = \la \phi(\bfx), \bfa \ra
$ where $\la \cdot, \cdot \ra$ is the Lie inner product.  Hence
\[ f(\bfx) = \la \phi(\bfx),
  \bfa \ra = L_\bfa(\phi(\bfx)) = (L_\bfa\circ \phi ) (\bfx)\]
where $L_\bfa$ is the linear function from $\boldR^{d+3}$ to $\boldR$
with $L_\bfa(\bfy) = \la \bfy, \bfa \ra $.  This we have shown $f =
L_\bfa\circ \phi.$

Conversely, given a linear function $L_\bfa: \boldR^{d+3} \to \boldR$
with $L_\bfa(\bfy) = \la \bfy, \bfa \ra $, we may partially invert
Equation \eqref{eqn: linear cov}, ignoring the $a_{d+3}$ term, to
associate to $\bfa$ the function $f(\bfx)$ as in Equation \eqref{eqn:
f} which has the property $f = L_\bfa\circ \phi.$
\end{proof}

\subsection{Well-known diagrams obtained from  Theorem \ref{thm: affine on quadric}  }

Many kinds of diagrams meet the hypotheses of Theorem \ref{thm:
  affine on quadric}, including some familiar ones. 
\begin{prop}\label{prop: affine diagrams}
  The following kinds of diagrams can be obtained as minimization
 diagrams as described in Theorem \ref{thm: affine on quadric}:
\begin{enumerate}
\item  The classical Voronoi diagram for point sites in $\boldR^d$;
\item Farthest point Voronoi diagrams for point sites in $\boldR^d$;
\item Power diagrams defined by a set of spheres in $\boldR^d$; 
\item Multiplicatively weighted Voronoi diagrams for point sites in
  $\boldR^d$; and
  \item The medial axis of a convex polygon. 
  \end{enumerate}
\end{prop}

\begin{proof}
Let $\calP = \{\bfp_1, \bfp_2, \ldots, \bfp_m\}$ be a set of points in
$\boldR^d.$ For each $i,$ let the function $f_i(\bfx) = \| \bfx -
\bfp_i \|^2 $ be the squared distance to the point $\bfp_i$.

The classical Voronoi diagram with sites in $\calP$ is the 
minimization diagram for the family $\calF = \{f_i\}$ of squared
distance functions. By Proposition \ref{prop: of form}, these functions meet
the hypotheses of Theorem \ref{thm: affine on quadric}.

The farthest point Voronoi diagram with sites in $\calP$ is the
minimization diagram for $-\calF = \{-f_i \, : \, f_i \in \calF\}$,
the family of the negatives of the squared distance functions.

The power diagram for a set of spheres $\{S_i\}$ in $\boldR^d$, by
definition, is the minimization diagram for the power functions as
defined in Equation \eqref{eqn: power def}.  The form of the power
function $\Pow(\bfx,S) = \| \bfx - \bfq \|^2 - r^2$ for a sphere
$S=S(\bfq,r)$ meets the hypotheses of Theorem \ref{thm: affine on
quadric}.

The multiplicatively weighted Voronoi diagram for the points $\calP$
with positive weights $w_1, \ldots, w_m$ is the minimization diagram
for the nonnegative functions $\frac{1}{w_i} \| \bfx - \bfp_i\|$. (See
\cite{aurenhammer1987power}, Section 6.3).  Equivalently, it is
the minimization diagram for the functions $\frac{1}{w_i^2} \| \bfx -
\bfp_i\|^2$.  The squared functions meet the hypotheses of Theorem
\ref{thm: affine on quadric}.

Let $P$ be a convex polygon defined by inequalities $\bfx \cdot \bfa_i
\ge h$ for $i=1, \ldots, m$.  The medial axis is the minimization
diagram for the functions $\bfx \cdot \frac{1}{\| \bfa_i\|}\bfa_i - \frac{1}{\| \bfa_i\|}
h$.  These are of the form in the  hypotheses of Theorem \ref{thm: affine on quadric}.
\end{proof}

\subsection{Order $k$ minimization diagrams.}
In this subsection we will show that Theorem \ref{thm: affine on
  quadric}
applies to order $k$  minimization diagrams.  
Before we prove Proposition \ref{prop: order k}, we make some definitions and
observations involving order $k$ minimization diagrams.  Let $\calF =
\{f_1, f_2, \ldots, f_m\}$ be a set of continuous functions mapping
$\boldR^d$ to $\boldR$.  Let $k \in [m]$.  For any subset $\{i_1, i_2,
\ldots, i_k\}$ of $[m]$ of cardinality $k$, define the function
$f_{i_1, i_2, \ldots, i_k}: \boldR^d \to \boldR$ by
   \begin{equation}\label{eqn: f_i} f_{i_1, i_2, \ldots, i_k}(\bfx) = f_{i_1}(\bfx) + f_{i_2}(\bfx)
     + \cdots + f_{i_k}(\bfx) . \end{equation}
Let $\calF^k$ be the collection of the $(\begin{smallmatrix} m \\ k
\end{smallmatrix} )$ functions $f_{i_1, i_2, \ldots, i_k}$:
   \[\calF^k = \{f_{i_1,
  i_2, \ldots, i_k} \, : \,  1 \le i_1 < i_2 < \cdots i_k \le m \}.\]
The next proposition states that the cells of the order $k$ Voronoi diagram
$\Vor^k(\calF)$ for $\calF$ are the same as the cells in the standard minimization
diagram $\MD(\calF^k)$ for $\calF^k$.
\begin{prop}\label{prop: order k bijection}
Let $\calF = \{f_1, f_2, \ldots, f_m\}$ be a set of continuous
functions mapping $\boldR^d$ to $\boldR$, let $k \in [m]$, and let
$\calF^k $ be as defined above.  There is a one-to-one correspondence
between the cells in the order $k$ minimization diagram
$\Vor^k(\calF)$ for $\calF$ and the cells in the minimization diagram
$\MD(\calF^k)$ for $\calF^k$.  The cell for the $k$ functions
$f_{i_1}, f_{i_2}, \cdots , f_{i_k}$ in the order $k$ minimization
diagram for $\calF$ is equal to the minimization region for the single
function $f_{i_1, i_2, \ldots, i_k}$ in $\MD(\calF^k).$
 \end{prop}

 \begin{proof} Fix $\bfx\in \boldR^d$ and $k \in [m]$.  Suppose that
$\bfx$ is in the order $k$ Voronoi cell for $f_{j_1}, \ldots, f_{j_k}$ in
the order $k$ minimization diagram for $\calF,$ where $f_{j_1}(\bfx)
\le f_{j_2}(\bfx) \le \cdots \le f_{j_k}(\bfx).$ Extend this ordering
to include all the values $f_1(\bfx), \ldots, f_m(\bfx)$:
\[ f_{j_1}(\bfx) \le f_{j_2}(\bfx) \le \cdots \le f_{j_m}(\bfx),\]
where $(j_1, j_2, \ldots, j_m)$ is a permutation of the set $[m].$
 It follows directly from this that $f_{j_1, j_2, \ldots, j_k}(\bfx)
\le f_{i_1, i_2, \ldots, i_k}(\bfx) $ for all subsets $\{i_1, i_2,
\ldots, i_k\}$ of $[m]$.  Hence $\bfx$ is in the minimization region for $f_{j_1, j_2, \ldots, j_k}(\bfx) $.

For the converse, suppose that $\bfx$  minimizes
for the function $f_{j_1, j_2, \ldots, j_k}$.  Then
\begin{equation}\label{eqn: ineq} f_{j_1, j_2, \ldots, j_k}(\bfx) \le f_{i_1, i_2, \ldots, i_k}
(\bfx) \end{equation} for all subsets $\{i_1, i_2, \ldots, i_k\}$ of cardinality
$k$.  We want to show that $f_{j_1}(\bfx), f_{j_2}(\bfx), \ldots,
f_{j_k}(\bfx)$ are the $k$ smallest values of $f_i (\bfx),$ for $i \in
[m].$  Were this not true, then
there would be some $r \not \in \{j_1, \ldots, j_k\}$ so that
$f_r(\bfx) < f_{j_i}(\bfx)$ for some $i \in \{1, 2, \ldots, k\}.$
Without loss of generality, suppose that $f_r(\bfx) <
f_{j_1}(\bfx)$. But then $f_{r, j_2, \ldots, j_k}(\bfx) < f_{j_1, j_2,
\ldots, j_k}(\bfx)$, a contradiction to \eqref{eqn: ineq}.
 \end{proof}

Proposition \ref{prop: order k} is an immediate
consequence of Proposition \ref{prop: order k bijection}.

\begin{proof}[Proof of Proposition \ref{prop: order k}.]
   Let $\calF = \{f_i \}$ be a family of continuous real-valued
   functions of  form $a \| \bfx - \bfq \|^2 +
   \bfb \cdot \bfx + c$  as in the statement of the proposition.   Let $\calF^k =
   \{ f_{{i_1, i_2, \ldots, i_k}}\}$
   be the family of functions defined by $\calF$  as in Equation
   \eqref{eqn: f_i}.
   A function of  form $a \| \bfx - \bfq \|^2 +
   \bfb \cdot \bfx + c$ is a multinomial of degree two or less  whose
   order two terms are  $a x_1^2, \ldots, a x_d^2$ for some real
   number $a$.   The sum of functions of this form are
   still of the same form.  Hence the functions $f_{{i_1, i_2, \ldots,
       i_k}}$ are of that form, and Theorem \ref{thm: affine on
     quadric} applies to the minimization diagram $MD(\calF^k)$ for $\calF^k.$  By  
Proposition   \ref{prop: order k bijection}, the cells of the order $k$ Voronoi
diagram for $\calF$ are the cells of $MD(\calF^k).$
\end{proof}
\section{Algorithm}\label{section: algorithm}

The algorithm is for computing the minimization diagram for a set $S$ of  functions
given by linear conditions in Lie sphere coordinates for Lie spheres
in $\boldR^d.$
Following is the algorithm to compute the generalized Voronoi diagram.

\begin{algorithm}\label{fig: algorithm}
  \KwData{Sites $S$}
  \KwResult{Voronoi vertices and connectivity}
  Map $S$ to canonical Lie sphere coordinates\;
  Do change of variables on $S$\;
  Compute and map bounding box\;
  Compute and map feasible point using bounding box\;
  Compute convex hull of Lie sphere coordinates by intersecting half-spaces in $\boldR^{d+2}$\;
  Get 1-dimensional faces $f$ between convex hull vertices\;
  Filter out edges with two bounding box equations\;
  Solve linear inequalities to find the voronoi edges\;
  Intersect edges with Lie quadric to find voronoi vertices (not medial axis)\;
\end{algorithm}

Let $n$ be $|S|$. Lines 1 and 2 of the algorithm are $O(nd)$. Lines 3
and 4 are $O(d)$. The Lie sphere coordinates are in
$\boldR^{d+2}$. Line 5, the intersection of half-spaces, uses the
incremental convex hull algorithm, which is $\mathcal{O}(n\log{n} +
n^{\lfloor\frac{d+3}{2}\rfloor})$
(\cite[p. 59]{boissonnat2018geometric}). For Line 6, assume the sites
are in general position. Each vertex $p$ is ``included in exactly
$d+2$ of the bounding hyperplanes'' (\cite[p. 54]{boissonnat2018geometric}). Then the number of 1-dimensional faces incident to each vertex $p$ is $d+2$. If there are $n$ hyperplanes, then there are $n(d+2)/2 = \calO(nd)$ 1-dimensional faces. Line 7 is $O(nd)$. In Line 8, to find the Voronoi edges, we solve a system of $d+1$ linear inequalities of half-spaces for each 1-dimensional face which is $O(nd^2)$. Line 9 is $O(nd)$.
The algorithm is dominated by lines 5 and 8, yielding $\calO(n\log{n} +
n^{\lfloor\frac{d+3}{2}\rfloor} + nd^2)$.

\bibliographystyle{amsalpha}
\bibliography{paper}

\end{document}